\documentclass[10pt,a4paper,reqno]{amsart}
\usepackage{latexsym,amssymb,amsmath,amsthm,amsfonts,enumerate,verbatim,xspace,exscale}

\input xy
\xyoption{all}
\CompileMatrices
\UseComputerModernTips

\theoremstyle{plain}
\newtheorem{theorem}{Theorem}[section]
\newtheorem{lemma}[theorem]{Lemma}
\newtheorem{proposition}[theorem]{Proposition}

\theoremstyle{definition}

\newtheorem*{name}{\theoremname}
\newcommand{\theoremname}{testing}

\newtheorem{rema}[theorem]{Remark}

\newenvironment{rem}
    {\begin{rema}}
    {\end{rema}}

\newtheorem{remarka}[theorem]{Remark}

\newenvironment{remark}
    {\begin{remarka}}
    {\end{remarka}}

\newtheorem{exx}[theorem]{Example}%[theorem]

\allowdisplaybreaks[3]

\newcommand{\caps}[1]{\textup{\textsc{#1}}}

\providecommand{\bysame}{\makebox[3em]{\hrulefill}\thinspace}

\newcommand{\iso}{isomorphism\xspace}

\newcommand{\bb}{\mathbb}
\newcommand{\ov}[1]{\mbox{$\overline{#1}$}}
\newcommand{\up}{\upshape}

\newcommand{\x}{$\hfill\Box$}

\newcommand{\longto}{\longrightarrow}
\newcommand{\hookto}{\hookrightarrow}
\newcommand{\toto}{\twoheadrightarrow}
\def\vv<#1>{\langle#1\rangle}

\newcommand{\id}{\mbox{$\text{\up{id}}\,$}}
\newcommand{\pr}{\mbox{$\text{\up{pr}}$}}

\providecommand{\det}{\mbox{$\text{\up{det}}\,$}}

\newcommand{\dd}[2]{\mbox{$\frac{\partial #2}{\partial #1}$}}

\newcommand{\om}{\omega}
\newcommand{\Om}{\Omega}

\newcommand{\lam}{\lambda}
\newcommand{\Lam}{\Lambda}

\newcommand{\wt}[1]{\mbox{$\widetilde{#1}$}}

\newcommand{\bsc}{\mbox{$\bigsqcup$}}

\newcommand{\R}{\mbox{$\bb{R}$}}
\newcommand{\C}{\mbox{$\bb{C}$}}

\newcommand{\by}[2]{\mbox{$\frac{#1}{#2}$}}
\newcommand{\cinf}{\mbox{$C^{\infty}$}}
\providecommand{\set}[1]{\mbox{$\{#1\}$}}
\newcommand{\subeq}{\subseteq}

\newcommand{\curv}{\mbox{$\textup{Curv}$}}

\newcommand{\ver}{\mbox{$\textup{Ver}$}}
\newcommand{\hor}{\mbox{$\textup{Hor}$}}
\newcommand{\ann}{\mbox{$\textup{Ann}\,$}}
\newcommand{\fix}[1]{\mbox{$\textup{Fix}(#1)$}}

\newcommand{\spr}[1]{\mbox{$/\negmedspace/_{#1}$}}
\newcommand{\sporb}{\mbox{$/\negmedspace/_{\mathcal{O}}$}}
\newcommand{\momap}{momentum map\xspace}

\newcommand{\lie}[1]{\mbox{$\text{\up{Lie}}(#1)$}}

\newcommand{\gu}{\mathfrak{g}}
\newcommand{\ho}{\mathfrak{h}}

\newcommand{\no}{\mathfrak{n}}

\newcommand{\wo}{\mathfrak{w}}

\newcommand{\Ad}{\mbox{$\text{\upshape{Ad}}$}}

\newcommand{\orb}{\mbox{$\mathcal{O}$}}

\newcommand{\SO}{\mbox{$\textup{SO}$}}
\newcommand{\so}{\mbox{$\mathfrak{so}$}}

\newcommand{\WW}{\mbox{$\mathcal{W}$}}

\newcommand{\lorb}[1]{\mbox{${#1}_{(L_{0})}^{\mathcal{O}}$}}

\newcommand{\ine}{\mbox{$\mathbb{I}$}}

\newcommand{\aklm}{Alekseevsky, Kriegl, Losik, Michor\xspace}

\title[Spinning particles in a Yang-Mills field]
      {Spinning particles in a Yang-Mills field}
\author{Simon Hochgerner}
\address{Fakult\"at f\"{u}r Mathematik\\
  Universit\"{a}t Wien\\
  Nordbergstrasse~15\\
  A-1090 Vienna, Austria}
\email{simon.hochgerner@univie.ac.at} 
\urladdr{http://www.mat.univie.ac.at/$\thicksim$simon}
\thanks{This work is supported by Fonds zur F\"{o}rderung der
       wissenschaftlichen Forschung (FWF), Projekt P 17108-N04}
\date{February 3, 2006}

\begin{document}

\begin{abstract}
Suppose that a Lie group $G$ acts properly on a configuration
manifold $Q$. We study the symplectic quotient of $T^*Q$ with respect
to the cotangent bundle lifted $G$-action at an arbitrary coadjoint
orbit level $\orb$. In particular, if $Q=Q_{(H)}$ is of single orbit
type we show that the symplectic quotient of $T^*Q$ at $\orb$ can be constructed
through a minimal coupling procedure involving the smaller cotangent
bundle $T^*Q_H$, the symplectic quotient of $\orb$ at $0$ with respect to
the $H$-action, and the diagonal Hamiltonian $N(H)/H$-action
on these symplectic spaces. A prescribed connection on $Q_H\to
Q_H/(N(H)/H)$ then yields a computationally 
effective way of explicitly realizing the
symplectic structure on each stratum of the symplectic quotient of
$T^*Q$.  
In an example this result is
combined with the projection method to produce a 
stratified Hamiltonian system
with very well hidden symmetries.
\end{abstract}
\maketitle

%\thispagestyle{empty}

%\pagebreak
%\thispagestyle{empty}
%\tableofcontents
%\pagebreak

\section{Introduction}
\subsection{A brief history of the orbit bundle picture in
  mechanics}
A general discussion of the history of symplectic reduction and its
variants can be found in the overview article by Marsden and
Weinstein~\cite{MW01}. In particular this article contains historical
remarks on the bundle picture in mechanics. The essence of this bundle
picture is the following. Let $Q$ be a smooth configuration manifold
acted upon in a proper and free fashion by a Lie group $G$. Then this
action can be cotangent lifted to give a Hamiltonian $G$-action on
$T^*Q$ with \momap $\mu: T^*Q\to\gu^*$. If $\orb$ is a coadjoint orbit
in the image of $\mu$ then the symplectic quotient $\mu^{-1}(\orb)/G
=: T^*Q\sporb G$ is a smooth symplectic manifold. The orbit bundle
picture is the observation that, under the additional assumption of a
prescribed principal bundle connection (such as the mechanical
connection) on $Q\to Q/G$, 
one has a smooth
symplectic fiber bundle 
\[
 \orb\hookto T^*Q\sporb G\longto T^*(Q/G)
\]
as well as a description of the reduced symplectic form on the total
space of this bundle. (\cite{AM78,Kum81,MMR84,Mon84,MP00}) 
If the orbit consists of a single
point only, i.e., $\orb=\set{\lam}$, this implies that
there is a symplectomorphism 
$T^*Q\spr{\lam}G\cong T^*(Q/G)$ where $T^*(Q/G)$ is equipped with a
magnetic symplectic form. In particular, one thus recovers the Abelian
version of cotangent bundle reduction as developed earlier in
\cite{Sma70,Sat77}.

The bundle picture was motivated by the articles of
Sternberg~\cite{S77} and Weinstein~\cite{W78} on minimal coupling
and Yang-Mills potentials. These articles show how to obtain the
equations of motions for a particle in a Yang-Mills field in a
symplectic framework.

If the $G$-action on $Q$ is assumed to be proper but not necessarily
free the symplectic quotient $T^*Q\sporb G$
is not a smooth manifold but a stratified
symplectic space in the sense of \cite{SL91,BL97,OR04}.
The case of non-free $G$-actions has first been considered in
Montgomery~\cite{Mon83}. In this paper Montgomery uses the point
reduction scheme and provides certain conditions under which the point reduced
cotangent bundle $T^*Q\spr{\lam}G := \mu^{-1}(\lam)/G_{\lam}$ embeds
as a subbundle in $T^*(Q^{\lam}/G_{\lam})$. Hereby $Q^{\lam} =
\set{q\in Q: \vv<\lam,\gu_q>=0}$, and the conditions alluded to are
related to the single orbit type assumption of the present paper. 
Assuming the $G$-action on $Q$ to be infinitesimally free, Emmrich and
R\"omer~\cite{ER90}
 have found the symplectic quotient of $T^*Q$ to be an
orbifold.
More recently,
Schmah \cite{Sch04} has proved a cotangent bundle specific slice
theorem at elements $(q,p)\in T^*Q$ whose momentum value
$\mu(q,p)=\lam$ is fully isotropic, i.e., fixed by the
$\Ad^*(G)$-action. Under the same assumption on the momentum value
$\lam$, Perlmutter, Rodriguez-Olmos, and Sousa-Diaz~\cite{PRS03} were able to
describe the geometry of the stratification of the reduced phase space
$T^*Q\spr{\lam}G$. 
Requiring $Q$ to be of single orbit type (that is $Q=Q_{(H)}$), 
the author found in \cite{H04} a generalization of the bundle picture for general
momentum values. This point of view was also pursued in \cite{hr05} to
describe the singular Poisson reduced space $(T^*Q)/G$.

\subsection{Statement of results}
Let $G$ be a Lie group acting properly on a configuration manifold $Q$.
Suppose $Q=Q_{(H)}$ is of single orbit type, that is $G_q$ is
conjugate to $H$ within $G$ for all $q\in Q$. Consider the cotangent
bundle lifted action by $G$ on $T^*Q$. This action is proper by
assumption and Hamiltonian with \momap $\mu$. Let $\orb$
denote a coadjoint orbit in the image of $\mu$. Clearly the
$\Ad^*(H)$-action on $\orb$ is Hamiltonian, 
albeit in general not free,
and we shall denote the
singular symplectic quotient of $\orb$ by this action at the $0$-level by
$\orb\spr{0}H$. Further, let $Q_H=\set{q\in Q: G_q=H}$ be the symmetry
type submanifold of $Q$ on which there is an induced action by
$W:=N(H)/H$ which is free by construction.  Thus there is a
diagonal action by $W$ on the product $T^*Q_H\times\orb\spr{0}H$ which
is free, proper, and Hamiltonian. We show in Theorem~\ref{thm:red}
that $\mu^{-1}(\orb)/G = T^*Q\sporb G$ is isomorphic as 
a singular symplectic space to the symplectic quotient
\[
 (T^*Q_H\times\orb\spr{0}H)\spr{0}W
\]
of 
$T^*Q_H\times\orb\spr{0}H$ with respect to the $W$-action at the
$0$-level.
Moreover, it is shown that the smooth symplectic strata of $T^*Q\sporb
G$ can be computed in similar manner involving $T^*Q_H$, the smooth
symplectic strata of $\orb\spr{0}H$, and the induced $W$-action. If
there is a principal bundle connection $\mathcal{A}$ given on $W\hookto Q_H\to
Q_H/W =: B$ 
(e.g., the mechanical connection)
then this minimal coupling 
construction 
yields a singular symplectic fiber bundle
\[
 \orb\spr{0}H
 \hookto
 (T^*Q_H\times\orb\spr{0}H)\spr{0}W
 \cong_{\mathcal{A}}
 (Q_H\times_B T^*B)\times_W\orb\spr{0}H
 \longto T^*B
\]
and, furthermore,
provides an explicit description of the
reduced symplectic structure on each of the symplectic strata of
$T^*Q\sporb G$. This exposition of the reduced symplectic structure 
is isomorphic to that given
in \cite{H04}. However, the advantage of the minimal coupling approach
is that it is computationally much more effective. For instance, in the
example presented in Section~\ref{sec:hidsym} it is virtually impossible
to obtain the reduced symplectic form by applying the bundle picture
construction of \cite{H04}. Another point of this proposed
construction is that it is specific to singular cotangent bundle
reduction. Indeed, if the $G$-action on $Q$ is free then the result
reduces to the shifting trick.  

The physical interpretation of the reduced space
$(T^*Q_H\times\orb\spr{0}H)\spr{0}W$ is that it is the phase space of
particles moving on $B$ in the presence of a Yang-Mills field and
subject to additional internal variables corresponding to
$\orb\spr{0}H$. 
These internal variables could, for example, be spin
variables. This is the case in the Calogero-Moser models considered in
\cite{AKLM03,H04,FP05}. Thus there is a qualitative difference between
singular ($H\neq\set{e}$) and regular ($H=\set{e}$) cotangent bundle
reduction. In the regular case the gauge group is $G$ and the internal
variables correspond to an $\Ad^{*}(G)$-orbit, while in the singular
case the relation between gauge group $W$ and internal variables
$\orb\spr{0}H$ is more intricate.

\begin{comment}
**************************************************************************
SO EIN BLOEDSINN:
As in the case of regular reduction, one can write down a set of Wong
equations (\cite{Won70}) 
on $(T^*Q_H\times\orb\spr{0}H)\spr{0}W\cong
T^*B\times_B(Q_H\times_W\orb\spr{0}H)$ 
and attempt to relate these to equations obtained
via Hamiltonian reduction of $T^*Q$. This is considered in
Remark~\ref{rem:wong}, and it thus transpires there is another
difference to regular cotangent bundle reduction: For the regular case
Montgomery~\cite{Mon84} has shown that the Sternberg~\cite{S77} and
Weinstein~\cite{W78} approach are in a sense equivalent to the
Kaluza-Klein approach of Kerner~\cite{Ker68} which is in turn shown to be
equivalent to Wong's equations. In the singular setting this is in
general not true. Indeed, following Theorem~\ref{thm:red} it will be
easy to see that the reduced function on $T^*Q\sporb G$ obtained from 
the Kaluza-Klein Hamiltonian on $T^*Q$
differs in general from the minimal coupling Hamiltonian of Sternberg
on $T^*Q\sporb G \cong (T^*Q_H\times\orb\spr{0}H)\spr{0}W$
by more than a Casimir function. 

Another remark (\ref{rem:quant}) is concerned
with possible applicability of Theorem~\ref{thm:red} to geometric
quantization of the singular symplectic quotient $T^*Q\sporb G$. 
***************************************************************************
\end{comment}

In Section~\ref{sec:strata} we make some comments on the general
problem of cotangent bundle reduction. This is the problem of
understanding the symplectic quotient $T^*Q\sporb G$ for general
proper $G$-actions on the configuration space $Q$. The view of
Section~\ref{sec:strata} is, in principle, that one should first
compute the spaces $T^*(Q_{(H)})\sporb G$ via Theorem~\ref{thm:red}
and then proceed by a case to case study to obtain the full symplectic
quotient $T^*Q\sporb G$. 

In Section~\ref{sec:hidsym} these ideas are applied to the diagonal
action of $G=\SO(5)$ on $Q=S^9\subset\R^5\times\R^5$. The
stratification of $(T^*S^9)/\SO(5)$ is exhibited according to the
remarks in Section~\ref{sec:strata}. Then a special coadjoint orbit
level $\orb\subset\so(5)^*$ is fixed. Employing Theorem~\ref{thm:red}
the reduced space $T^*S^9\sporb\SO(5)$ together with the induced
symplectic form on each of its strata is computed. Moreover, this
reduction process is carried out in the presence of a $\SO(5)$-invariant
Hamiltonian function on $T^*S^9$ whence we arrive at a stratified
Hamiltonian system. In fact, the example is chosen so that there are
only two strata and the induced Hamiltonian system on each of these is
described shortly. The system on the small stratum corresponds to the
motion of a charged particle on a closed disk under the influence of
an electromagnetic field. The system on the big stratum is more
complicated and a physical interpretation is attempted at the end of
Section~\ref{sec:hidsym}.

%\subsection{Relation with other work}
%Possible remarks about implications for quantization.

\section{Preliminaries and notation}
All manifolds to be considered are Hausdorff, paracompact, 
finite dimensional, and smooth in the $\cinf$-sense.
Let $(M,\om)$ be a Hamiltonian $G$-space, i.e.\ $(M,\om)$ is a
symplectic manifold acted upon via symplectomorphisms 
by a Lie group
$G$ such that there is an equivariant \momap $J: M\to\gu^*$. Whenever
the \momap is clear from the context we will write $M\spr{\lambda}G$
for $J^{-1}(\lam)/G_{\lam}$ and $M\sporb G$ for $J^{-1}(\orb)/G$ where
$\orb$ is a coadjoint orbit. We will only be concerned with left
actions by compact Lie groups. 
If the action is written as $l: G\times M\to M$, $(k,x)\mapsto
l(k,x)=l_{k}(x)=l^{x}(k)=k.x$ 
we can tangent bundle lift it via
$k.(x,v):=(k.x,k.v):=Tl_{k}.(x,v)=(l_{k}(x),T_{x}l_{k}.v)$ 
for $(x,v)\in TM$ to an action on $TM$. 
As the
action consists of transformations by diffeomorphisms 
it may also be lifted to
the cotangent bundle. This is the cotangent lifted action which is
defined by
$k.(x,p):=(k.x,k.p):=T^{*}l_{k}.(x,p)=
(k.x,T_{k.x}^{*}l_{k^{-1}}.p)$ where $(x,p)\in T^{*}M$. 
Our notation for the fundamental vector field is 
$
 \zeta_{X}(x) 
 := \dd{t}{}|_{0}l(\exp(+tX),x)
  = T_{e}l^{x}(X)
$
where $X\in\gu$.

If the action by $G$ on $M$ is proper, in the sense that $G\times M\to
M\times M$, $(k,x)\mapsto(x,l(k,x))$ is a proper mapping, then we have
the Slice and Tube Theorem at our disposal. 
An exposition of these
facts can be found in Palais and Terng~\cite{PT88}, 
Duistermaat and Kolk~\cite{DK99}, or
Kawakubo~\cite{kaw91}, 
for example. We say $M$ is a proper $G$-manifold if $G$ is a Lie
group acting properly on $M$. 

Let $x\in M$ and $H:=G_x$. Then we define 
$(H) := \set{kHk^{-1}: k\in G}$ to be the $G$-\caps{conjugacy class}
of $H$, and say that $x$ is of \caps{isotropy type} $(H)$. 
If $G$ is compact we can define a 
a partial ordering on the set of all  
conjugacy classes (i.e., isotropy types) of the $G$-action on $M$ as follows.
Namely, say $(H)\prec (L)$
if $L$ is conjugate to a subgroup of $H$.
 
The submanifold of all points $x\in M$ such that $(G_x)=(H)$ is
denoted by $M_{(H)}$, 
that is,
\[
 M_{(H)}
 =
 \set{q\in M: G_q \text{ is conjugate to } H \text{ within } G},
\]
and 
this is called the \caps{isotropy type} or \caps{orbit type}
submanifold of $M$ of type $(H)$.
Further, let
\[
 M_H
 := 
 \set{q\in Q: G_q=H}
\text{ and }
 M^H
 := 
 \set{q\in Q: H\subset G_q}
\]
denote the \caps{symmetry type} and \caps{fixed point} submanifold,
respectively, of $H$.

\begin{theorem}\label{thm:open-dense}
Suppose $G$ is compact and $M$ is a connected $G$-manifold. Then the
following are true.
\begin{enumerate}[\up (1)]
\item
There exists a unique maximal isotropy type $(H)$ characterized by
requiring that
$\dim H = \inf\set{\dim G_x: x\in M}$ and that $H$ has the least
number of connected components among all those isotropy subgroups of $G$
which have dimension equal to $\inf\set{\dim G_x: x\in M}$.
\item
If $(H)$ is the maximal isotropy type then $M_{(H)}$ is open and dense
in $M$ and $M_{(H)}/G$ is connected. 
\end{enumerate}
\end{theorem}

\begin{proof}
The proof is based on the Palais' Slice Theorem and can be found in,
e.g., Palais and Terng \cite{PT88}. 
\end{proof}

Throughout the paper we 
let $Q$ be a connected manifold, and $G$ a connected Lie group which
acts properly from the left on $Q$. 
We will write the action alternatively as $l(k,q)=l^q(k)=l_k(q)=k.q$
where $k\in G$ and $q\in Q$.
We will further assume throughout that 
$Q_{(H)}$
is non-empty where $(H)$ is a fixed isotropy type of the $G$-action. 
Clearly
$G.Q_H=Q_{(H)}$, and we have an induced action by $G$ on $Q_{(H)}$
which is of single isotropy type, and we have an induced free action
by $W=W(H):= N(H)/H$ on $Q_H$. Notice, moreover, that right
multiplication in the group induces a left action by $W$ on $G/H$,
that is $w.[k]=[kw^{-1}]$. The following description of $Q_{(H)}$ 
is the key ingredient of the sequel. 

\begin{theorem}[Structure Theorem]\label{thm:str}
The orbit projection $Q_{(H)}\toto Q_{(H)}/G$ is a smooth fiber bundle
with typical fiber $G/H$. Moreover, it is diffeomorphic to
the bundle associated to the
principal bundle $Q_H\toto Q_H/W$ with respect to the left action by
$W$ on $G/H$. The diffeomorphism is given by $\kappa: Q_H\times_{W}G/H\to
Q_{(H)}$, $[(q,kH)]_W\mapsto k.q$.
\end{theorem}

\begin{proof}
See Duistermaat and Kolk \cite{DK99}. 
\end{proof}

\begin{rem}
In Palais and Terng~\cite{PT88} the above theorem is proved for proper
Fredholm Riemannian $G$-manifolds. 
These are Riemannian manifolds $M$ which are modeled on a Hilbert
space such that the $G$-action is proper and isometric, and such that
the tangent map at each point of the orbit projection mapping $M\toto
M/G$ is a Fredholm map of Hilbert spaces. It should thus be possible
to generalize Theorem~\ref{thm:red} to this setting.
\end{rem}

\begin{rem}\label{rem:conn}
Let $(M,\om)$ be a Hamiltonian $G$-space such that the $G$-action on
$M$ is proper and the \momap is denoted by $J$. 
By \cite{SL91,BL97,OR04} the symplectic quotient $J^{-1}(\orb)/G =
M\sporb G$ is a stratified symplectic space. This means that the
reduced space is a Whitney (B)-stratified space as defined in
\cite{Mat70,DK99}, 
its strata are smooth symplectic manifolds, and the
inclusion mapping of each stratum into $M\sporb G$ is a Poisson
morphism with respect to the function Poisson structure on $M\sporb
G$.
If $J$ is not equivariant it is for reasons explained in \cite{OR04}
important to consider the stratification of $M\sporb G$ by the
connected components of the smooth symplectic manifolds
$(J^{-1}(\orb)\cap M_{(L)})/G$ as opposed to considering the
stratification given by the disconnected pieces 
$(J^{-1}(\orb)\cap M_{(L)})/G$. Here $(L)$ denotes a conjugacy class of
some isotropy subgroup of the $G$-action. However, it is notationally
simpler to deal with the stratification given by the disconnected
pieces, and since all momentum map appearing throughout this paper
will be equivariant by construction we will regard the disconnected
symplectic manifolds $(J^{-1}(\orb)\cap M_{(L)})/G$ as the strata of the
symplectic quotient. One can, of course, at any point pass to the
finer stratification given by  the connected components of these
strata. 
\end{rem} 

On $T^*Q$ we have a canonical \momap $\mu: T^*Q\to\gu^*$ given by
$\vv<\mu(q,p),X> = \vv<p,\zeta_X(q)>$ where $\vv<.,.>$ denotes the
dual pairing in the appropriate sense, and this is referred to as the
cotangent bundle \momap of the lifted $G$-action. 
We will denote      
the cotangent bundle \momap on $T^*Q_{(H)}$ by the same symbol $\mu$,
and, moreover, the $W=N(H)/H$-\momap on $T^*Q_H$ will be denoted by
$\mu$ as well. This will not cause any confusion since the meaning
will be clear from the context.

\section{Description of the big phase space}\label{sec:gen}
The above Theorem \ref{thm:str} is expressed in the following diagram.
\[
\tag{D1}
\xymatrix{
 {Q_{(H)}}
  \ar @{<-}[r]^-{\kappa}_-{\cong}
  \ar @{->>}[d]
 &
 {Q_H\times_{W}G/H}
  \ar @{<<-}[r]
  \ar[d]
 & 
 {Q_H\times G/H}
  \ar[d]^-{\textup{pr}_1}\\
 {Q_{(H)}}\ar @{<-}[r]^-{\cong}
 &
 {Q_H/W}\ar @{<<-}[r]
 &
 {Q_H}
}
\]
The point to be exploited in the sequel is that $\kappa$ lifts to a
symplectomorphism of the cotangent bundles $T^*Q_{(H)}$ and
$T^*(Q_H\times_{W}G/H)$ which is equivariant with respect to the lifted
$G$-actions. 
Using commuting reduction (see Marsden et al.\ \cite{MMOPR03}) 
in order to two times 
interchange Poisson reduction via left multiplication with
two-shot symplectic reduction via right multiplication we thus get 
\begin{align*}
 (T^*Q_{(H)})/G
 &\cong
 T^*(Q_H\times_{W}G/H)/G\\
 &= 
 (T^*Q_H\times(T^*G)\spr{0}H)\spr{0}W/G
 =
 (T^*Q_H\times\ann\ho/H)\spr{0}W.
\end{align*}
Hereby, as above, $H$ and $W$ act on $G$ by inversion of right
multiplication, and these actions are cotangent lifted to $T^*G$. 
Thus the problem of understanding $(T^*Q_{(H)})/G$ reduces is to that
of understanding $(T^*Q_H\times\ann\ho/H)\spr{0}W$ which is much
easier since, firstly, $W$ acts on $T^*Q_H$ freely and, secondly, the
involved spaces are smaller. Moreover, in many examples $W$ is an
Abelian and sometimes even finite group. 
Notice also that we can equip $\gu$ with an $H$-invariant inner product,
and thus identify 
\[
 \wo=\lie{W} 
 \cong
 \lie{N(H)}/\ho
 \cong
 \lie{N(H)}\cap\ho^{\bot}
 =
 \ho^{\bot}\cap\fix{H}=(\ho^{\bot})^H.
\]
That is, $H$ does not act on
$\wo$.
We will assume this identification tacitly for the rest of the paper.

\begin{lemma}\label{lem:J_W}
The \momap $J_W:
T^*Q_H\times\ann\ho/H\to\wo^*$ of the $W$-action is given by 
\[
 J_W(q,p,[\lam])=\mu(q,p)-\lam|\wo
\]
where $\mu$ is the cotangent bundle \momap on $T^*Q_H$ with respect to
the $W$-action. 
\end{lemma}

\begin{proof}
Indeed, notice firstly that $\ann\ho/H = \gu^*\spr{0}H$ is the Poisson
reduced space of $\gu^*$ with respect to the Hamiltonian $H$-action at
$0\in\ho^*$. 
More precisely, the $H$-action is given by
$h.\lam=\Ad^*(h^{-1}).\lam=\lam\circ\Ad(h)$. 
Now, the Hamiltonian $N(H)$-action which is given by the same formula
on $\gu^*$ induces an Hamiltonian action on
$\ann\ho/H$. 
This induced action is the action by $W=N(H)/H$ under consideration.
Its \momap thus computes straightforwardly to be $\ann\ho/H\to\wo^*$,
$[\lam]\mapsto-\lam|\wo$.  
\end{proof}

Let $\ver$ denote the vertical subbundle of $TQ_H$ with respect to
the the principal bundle projection
$Q_H\toto B$, and $\hor^*$ shall denote the dual horizontal subspace,
i.e.\ the subspace of covectors which annihilate horizontal
vectors. 
Assume there is a principal bundle connection form $\mathcal{A}:
TQ_H\to\wo$ given. 
Using this connection form we define $\hor$ and
$\ver^*$ in the usual way. 
Since $\mathcal{A}$ reproduces generators of fundamental vectorfields
the definition of the cotangent bundle \momap $\mu: TQ_H\to\wo^*$
implies that
$\hor^*_q=\mu_q^{-1}(0)$,
and, further, 
$\mathcal{A}_q^*: \wo^*\to\ver^*_q$ provides an
inverse to $\mu_q|\ver^*_q$. Thus the zero level set of $J_W$ turns
out to be
\begin{align*}
 J_W^{-1}(0)
 &\cong
 \set{(q,p_0+\mathcal{A}_q^*(\lam_0),[\lam_0+\lam_1]_H):
      p_0\in\hor_q^*,\lam_0\in\wo^*,\lam_1\in\ann(\ho+\wo)}\\
 &\cong
 \hor^*\times\wo^*\times\ann(\ho+\wo)/H.
\end{align*}
Therefore, we get the following description of the reduced space
$(T^*Q_{(H)})/G$.

\begin{proposition}\label{prop:cot-red}
Under the above assumptions $(T^*Q_{(H)})/G$ is isomorphic as a
stratified space to
$(Q_H\times_{B}T^*B)\times_{W}(\wo^*\times\ann(\ho+\wo)/H)$ where
$B=Q_H/W\cong Q_{(H)}/G$. Moreover,
\[
 \ann\ho/H
 =
 \wo^*\times\ann(\ho+\wo)/H
 \hookto(Q_H\times_{B}T^*B)\times_{W}(\wo^*\times\ann(\ho+\wo)/H) 
 \to T^*B
\]
is a singular fiber bundle in the sense of \cite{hr05}.
\end{proposition}

\begin{remark}
We can use the \iso of the above  proposition to endow
$(Q_H\times_{B}T^*B)\times_{W}(\wo^*\times\ann(\ho+\wo)/H)$ with a
Poisson structure. The induced 
Poisson structure on $(T^*Q_{(H)})/G$ was
described in Hochgerner and Rainer \cite{hr05}. 
More explicitly,
notice that the $G$-equivariant diffeomorphism $\kappa:
Q_H\times_{W}G/H\to Q_{(H)}$ 
lifts to induce a $G$-equivariant diffeomorphism
$\kappa_{0}: 
 (Q_H\times G\times_H\ann\ho)/W\to\bsc_{q\in Q_{(H)}}\ann\gu_q
 = \ver^*(Q_{(H)}\toto Q_{(H)}/G)$. 
Thus we get an induced \iso
$Q_H\times_W\ann\ho/H\cong(\bsc_{q\in Q_{(H)}}\ann\gu_q)/G$ of
stratified spaces. Moreover, pulling back commutes with forming
associated bundles, that is
$(Q_H\times_{B}T^*B)\times_W(G\times_H\ann\ho) 
\cong
T^*B\times_B(Q_H\times_W(G\times_H\ann\ho))$.  
We therefore get
\begin{multline*}
 (Q_H\times_{B}T^*B)\times_{W}(\wo^*\times\ann(\ho+\wo)/H)
 \cong\\
 \cong
 T^*B\times_B(Q_H\times_W\ann\ho/H)
 \cong
 T^*(Q_{(H)}/G)\times_{Q_{(H)}/G}(\bsc_{q\in Q_{(H)}}\ann\gu_q)/G
\end{multline*} 
where the last space is the Weinstein space description of $(T^*Q_{(H)})/G$
whose Poisson structure is described in \cite[Theorem 5.11]{hr05}. 
The first \iso in the above equation can be seen as a singular
instance of the Weinstein picture being equivalent to the Sternberg
description. (See Perlmutter and Ratiu \cite{PR04} for the $\cinf$-regular
version of these descriptions.)
\end{remark}

\section{The reduced phase space}\label{sec:curv}
Let $\pi_W: Q_H\toto Q_H/W=B$ and $\rho: Q_H\times G/H\to
Q_H\times_{W}G/H$ denote the orbit projections.

\subsection{The mechanical connection}\label{sub:mech-con}
Assume $G$ acts by isometries on $Q$ with respect to a Riemannian
metric $g$. We denote the restriction of $g$ to the totally geodesic
submanifold $Q_H$ by $g$ again. Thus $W$ acts by isometries on
$(Q_H,g)$. 

In this setting $\mu: T^*Q_H\to\wo^*$ 
yields a natural connection form on the principal
bundle $W\hookto Q_H\toto Q_H/W=:B$ as follows. 
Define the
\caps{moment of inertia tensor} $\ine: Q_H\to\wo^*\otimes\wo^*$ by
$\ine_q(X,Y)=g_q(\zeta_X(q),\zeta_Y(q))$ where $X,Y\in\wo$. 
Since the $W$-action on $Q_H$ is
by isometries it is clear that
$\ine_w.q(\Ad(w).X,\Ad(w).Y)=\ine_q(X,Y)$ for all $w\in W$,
whence $\ine$ defines a smooth family of inner products on $\wo$
depending on $q\in Q_H$.  
The natural connection form to be thus constructed is the
\caps{mechanical connection} $\mathcal{A}$ defined by
\[
\tag{D2}\label{d:A}
\xymatrix{
 {T_qQ_H}\ar[r]^-{\mathcal{A}_q}\ar[d]^-{\cong}_-{g_q}
 & 
 {\wo}\\
 {T_q^*Q_H}\ar[r]^-{\mu_q}
 &
 {\wo^*}\ar[u]_-{\ine_q}^-{\cong}
}
\]
Let $\Phi :=
\zeta^W\circ\mathcal{A}\in\Om^1(Q_H;TQ_H)$ be the principal connection
associated to the connection form $\mathcal{A}$.

\subsection{The mechanical curvature}\label{sub:mech-curv}
Continue to assume that $G$ acts on $Q$ by isometries.
The following diagram shows that the generalized mechanical connection $A$
of the bundle $\pi: Q_{(H)}\toto Q_{(H)}/G$ as defined in \cite{H04}
is associated to the mechanical connection $\mathcal{A}$ of 
Subsection~\ref{sub:mech-con}. 
The term generalized connection is to be understood in the context of
Alekseevsky and Michor~\cite{AM95}.
\[
\xymatrix{
 {TQ_H\times T(G/H)}
  \ar[r]^-{\Phi\times\textup{id}}\ar[d]_-{T\rho}&
 {\ver(\pi_W)\times T(G/H)}
  \ar[r]^-{A\times\textup{id}}_-{\cong}\ar[d]^-{T\rho}&
 {Q_H\times\wo\times G\times_H\ho^{\bot}}
  \ar[d]^-{\rho_0}\\
 {T(Q_H\times_{W}G/H)}
  \ar[r]^-{\widetilde{\Phi}}\ar[d]_-{T\kappa}^-{\cong}&
 {\ver(Q_H\times_{W}G/H\to B)}
  \ar[r]^-{\widetilde{\mathcal{A}}}_-{\cong}&
 {\by{(Q_H\times\wo\times G\times_H(\ho+\wo)^{\bot})}{W}}
  \ar[d]^-{\cong}\\
 {TQ_{(H)}}
  \ar[r]^-{A}&
 {\bsc_{q\in Q_{(H)}}{\gu/\gu_q}}&
 {(Q_H\times G\times_H\ho^{\bot})/W}
  \ar[l]_-{\kappa_{0}}^-{\cong}
}
\]
where $\rho_0$ is orthogonal projection composed with orbit
projection.
In particular, by construction, 
$\Phi\times\id$ and 
$(\kappa^{-1})^*\wt{\Phi} = \zeta^G\circ A$ 
are $\kappa\circ\rho$-related. Thus the same is true for the
respective curvatures. That is, the curvatures 
\[
 R_{\Phi\times\textup{id}}
 = 
 \by{1}{2}[\Phi\times\id,\Phi\times\id]
 = 
 \by{1}{2}[\Phi,\Phi]\times0
 =
 R_{\Phi}\times0
 =
 -\zeta^W\circ\curv^{\mathcal{A}}\times0
\]
and
\[
 R_{\widetilde{\Phi}}
 =
 \by{1}{2}[\widetilde{\Phi},\widetilde{\Phi}]
 =
 -\zeta^G\circ\curv^{A}
\]
are $\kappa\circ\rho$-related. The bracket that appears here is the
Fr\"olicher-Nijenhuis bracket,
and relatedness means that
$T(\kappa\circ\rho)\circ(R_{\Phi}\times0)
 =
 R_{\widetilde{\Phi}}\circ(T(\kappa\circ\rho)\oplus  T(\kappa\circ\rho))$.
For the curvature forms this implies 
\[
 T(\kappa\circ\rho)\circ((\zeta^W\circ\curv^{\mathcal{A}})\times0)
 =
 \zeta^G\circ\curv^A\circ\oplus^2 T(\kappa\circ\rho).
\]
Thus we arrive at the following assertion which, for emphasis, we
record as a proposition.

\begin{proposition}\label{prop:curv}
The generalized mechanical curvature $\curv^A$ is given by the formula
\[
 \curv^A\circ\oplus^2 T(\kappa\circ\rho)
 = 
 \kappa_{0}\circ\rho_0\circ(\curv^{\mathcal{A}}\times0).
\]
Therefore, 
$\curv^A: \Lam^{2}TQ_{(H)}\to\gu$ is $G$-equivariant, and
$\Ad(h).\curv^A(v,w)(q) = \curv^A(v,w)(q)$ for all $h\in G_q$.
Moreover, $\curv^{\mathcal{A}}$ is horizontal and thus drops to a
closed two-form $\curv_0^{\mathcal{A}}$ 
on $B=Q_H/W$ with values in the adjoint bundle
$Q_H\times_W\wo$, and the isomorphisms $B\cong Q_{(H)}/G$,
$\smash{Q_H\times_W\wo \cong (\bsc_{q\in
  Q_{(H)}}\fix{G_q}\cap{\gu/\gu_q})/G}$ relate $\smash{\curv_0^{\mathcal{A}}}$
and $\smash{\curv_0^A}$. 
\end{proposition}

\begin{remark}
Notice that
this equation greatly facilitates the work one has to do in computing
the mechanical curvature $\curv^A$ in examples. Furthermore, it gives
a geometrically satisfactory explanation of the otherwise somewhat
surprising properties of \cite[Proposition 4.1]{hr05}.
\end{remark}

\begin{remark}
Of course, Proposition~ \ref{prop:curv} is valid for the curvature form
form $\curv^A$ of any connection form $A$ associated to any 
principal
(not
necessarily mechanical) connection form $\mathcal{A}$ on $\pi_W:
Q_H\toto B$. However, since in applications we are concerned with
mechanical connections only we chose to state it in this way.  
\end{remark}

\subsection{Cotangent bundle reduction via minimal coupling}
Combining Propositions \ref{prop:cot-red} and \ref{prop:curv} with the
results of \cite[Section 6]{hr05} one is lead to expect 
that the singular symplectic
leaves of the Poisson reduced phase space $(T^*Q_{(H)})/G$ are given
by spaces of the form  
$(Q_H\times_{B}T^*B)\times_W\orb\spr{0}H$ where $\orb$ is a coadjoint
orbit of $G$ and $\orb\spr{0}H = (\orb\cap\ann\ho)/H$. The adjective
singular means that these leaves are actually stratified symplectic
spaces, and the smooth symplectic leaves of $(T^*Q_{(H)})/G$ are thus
given by the connected components of the smooth symplectic strata of
$(Q_H\times_{B}T^*B)\times_W\orb\spr{0}H$. 

For the following, let
\[
 (\orb\spr{0}H)_{(L_0)^H}
 :=
 (\orb_{(L_0)^H}\cap\ann\ho)/H
\]
where $\orb_{(L_0)^H}$ denotes the set of $\lam\in\orb$ such that
$H_{\lam}=G_{\lam}\cap H$ is conjugate to $L_0\subset H$ within
$H$. By virtue of symplectic reduction the smooth manifold
$(\orb\spr{0}H)_{(L_0)^H}$ inherits a symplectic form $\lorb{\Om}$ from the
(positive) KKS-form on the coadjoint orbit $\orb$. 
Further, we define
\[
 (\orb\spr{0}H)_{(L_0)^{N(H)}}
 :=
 (\orb_{(L_0)^{N(H)}}\cap\ann\ho)/H
\]
where $(L_0)^{N(H)}$ denotes the conjugacy class of $L_0$ in $N(H)$.

\begin{lemma}\label{lem:strat-pres}
The following are equivalent.
\begin{enumerate}[\up (1)]
\item
For all $n\in N(H)$ there is an $h\in H$ such that $nL_0n^{-1} =
hL_0h^{-1}$, that is, $(L_0)^{N(H)}=(L_0)^H$.
\item
$N_{N(H)}(L_0)/N_H(L_0) = W$.
\end{enumerate}
If these conditions are satisfied then
the action by $N(H)$ on $\orb$ 
induces a Hamiltonian action by $W$ on $(\orb\spr{0}H)_{(L_0)^H}$.

There is always an induced action by 
$N_{N(H)}(L_0)/N_H(L_0)$ on $(\orb\spr{0}H)_{(L_0)^H}$, and this
action is Hamiltonian.
\end{lemma}

\begin{proof}
This follows from the fact that $H_{n.\lam} = nH_{\lam}n^{-1}\subset
H$ for all $n\in N(H)$ and $\lam\in\orb$.
\end{proof}

In view of the previous lemma we will need to consider the
space
$W.(\orb\spr{0}H)_{(L_0)^H}$ which we will refer to as the $W$-sweep of
$(\orb\spr{0}H)_{(L_0)^H}$.

\begin{lemma}\label{lem:W-sweep}
The $W$-sweep of $(\orb\spr{0}H)_{(L_0)^H}$ is a 
coproduct 
\[
 \bsc_{[n]\in W}(\orb\spr{0}H)_{(nL_0n^{-1})^H}
\]
where $[n]=nH\in N(H)/H=W$, 
and has the following properties.
\begin{enumerate}[\up (1)]
\item
It is a finite disjoint union of strata of 
$\orb\spr{0}H$ 
(possibly with finitely many connected components) 
all of which are symplectomorphic to 
$(\orb\spr{0}H)_{(L_0)^H}$.
\item
$W.(\orb\spr{0}H)_{(L_0)^H} = (\orb\spr{0}H)_{(L_0)^{N(H)}}$
\item
When $\orb$ is compact
then there is exactly one regular stratum in $\orb\spr{0}H$. It is open,
dense and connected, and, furthermore, preserved by the $W$ action.  
\end{enumerate}
\end{lemma}

\begin{proof}
Indeed, as in the proof of Lemma \ref{lem:J_W} we notice that $W$ acts
by Poisson morphisms on $\orb\spr{0}H$. It is a general fact (e.g.,
\cite{SL91}) that homeomorphisms of singular Poisson spaces which are 
Poisson are also strata preserving. Therefore, $W$ maps strata onto
strata,
and it is easy to verify that 
$w=[n]=nH$ maps
$(\orb\spr{0}H)_{(L_0)^H}$
symplectomorphically onto
$(\orb\spr{0}H)_{(nL_0n^{-1})^H}$. 
Therefore, the space under consideration is a coproduct of the
asserted form, and it is even 
finite since by compactness
of $H$ there are only
finitely many strata of $\orb\spr{0}H$. 
This proves (1).

(2)
This assertion is straightforward to verify.

(3)
By a result of Kirwan \cite{Kir84} the pre-image of $0$ of a proper
\momap is always connected. Since $\orb$ is compact it thus follows
that $\orb\cap\ann\ho$ is connected. Therefore, by \cite{SL91} the
reduced space has a unique open, dense and connected stratum which is
characterized by Theorem \ref{thm:open-dense}. These characterizing
properties are preserved by the $W$-action -- whence the assertion.
\end{proof}

The following theorem uses the notion of a Hamiltonian fiber bundle
which is defined in Subsection \ref{sub:min} below.

\begin{theorem}[Cotangent bundle reduction as minimal coupling]\label{thm:red}
Suppose $G$ is a connected Lie group acting properly on a connected
manifold $Q$, 
and let $(H)$ be an isotropy type of this action.  
Assume that $\orb$ is a coadjoint orbit 
contained in the image of the cotangent bundle \momap $\mu:
T^*Q_{(H)}\to\gu^*$. 
The following are true.
\begin{enumerate}[\up (1)]
\item
There is an \iso of singular symplectic spaces identifying
the symplectic reduced space 
\[
 (T^*Q_{(H)})\sporb G :=
 \mu^{-1}(\orb)\cap(T^*Q_{(H)})/G
\] 
and the symplectic reduced
space at $0\in\wo^*$ with respect to the free diagonal action by
$W=N(H)/H$ on $T^*Q_H\times(\orb\spr{0}H)$. 
\item
Let $\mathcal{A}$ be a principal connection form on $\pi_W\toto B$.
Then $\mathcal{A}$ yields a
singular fiber bundle 
\[
 \orb\spr{0}H\hookto
 (T^*Q_H\times\orb\spr{0}H)\spr{0}W\cong
 (T^*Q_{(H)})\sporb G\longto
 T^*B
\]
in the sense of \cite{hr05}.
The transition functions of this bundle take values in $W$ which acts
on $\orb\spr{0}H$ by Hamiltonian transformations. Thus the fiber bundle is
Hamiltonian.
\end{enumerate}
Suppose, further, that $(L)$ is an isotropy type of the lifted
$G$-action on $T^*Q_{(H)}$. Then $L$ is conjugate to a subgroup
$L_0\subset H$ within $G$ and the following hold.
\begin{enumerate}[\up (1)]\setcounter{enumi}{2}
\item\label{item:iso}
There is
a $\cinf$-symplectomorphism 
\[
 ((T^*Q_{(H)})\sporb G)_{(L)}
 \cong
 (T^*Q_H\times(\orb\spr{0}H)_{(L_0)^H})\spr{0}W
\]
where 
$T^*Q_H\times(\orb\spr{0}H)_{(L_0)^H}$ is equipped with the obvious
product symplectic form.
\item\label{item:sp}
Let $\mathcal{A}$ denote a principal connection form on $Q_H\toto
Q_H/W=B$ with curvature form $\curv^{\mathcal{A}}$.   
The symplectic form on $T^*Q_H\times W.(\orb\spr{0}H)_{(L_0)^H}$
restricted to the $0$ level set of the \momap $J_W$ of the diagonal
$W$-action
can be
described by the minimal coupling form
\[
 \pr_1^*\eta^*\Om^{B} 
  -
 \vv<\phi\circ\pr_2,(\eta\circ\tau\circ\pr_1)^{*}\curv_0^{\mathcal{A}}> 
  - \pr_2^*\lorb{\Om}
\]
where $\phi: W.(\orb\spr{0}H)_{(L_0)^H}\to\wo^*$ is the $W$-\momap
given by $\phi([\lam]_H) = \lam|\wo$.
Furthermore,
$\Om^B$ is the canonical symplectic form on $T^*B$, $\tau:
T^*B\to B$, 
$(\pr_1,\pr_2): T^*Q_H\times W.(\orb\spr{0}H)_{(L_0)^H}
\to T^*Q_H\times W.(\orb\spr{0}H)_{(L_0)^H}$ 
are the obvious projections,
$\eta: T^*Q_H\to T^*B$ is the projection defined by the connection
$\mathcal{A}\in\Om^1(Q_H;\wo)$,
and $\curv_0^{\mathcal{A}}$ is
the induced form on $B$ from the basic form $\curv^{\mathcal{A}}$. 
This coupling form is horizontal and drops to the induced 
symplectic form on the reduced space 
$J_W^{-1}(0)/W 
=
(T^*Q_H\times W.(\orb\spr{0}H)_{(L_0)^H})\spr{0}W$.
\item
Assume $G$ acts on $Q$ by isometries with respect to some Riemannian
structure. 
Then
the symplectomorphism of item \ref{item:iso}
is 
compatible with the Weinstein description of 
\cite{H04}
in the following sense.
Let $A$ be the generalized mechanical connection 
associated to $\mathcal{A}$ as in
Subsection \ref{sub:mech-curv}, and let $\WW=\WW(A)$ be the
$A$-dependent Weinstein realization of $T^*Q_{(H)}$ as in
\cite[Section 5]{H04}. Then the symplectic structure on 
$(T^*Q_H\times W.(\orb\spr{0}H)_{(L_0)^H})\spr{0}W$ 
of item \ref{item:sp}
and the induced symplectic form on
\begin{align*}
 (T^*Q_{(H)}\sporb &G)_{(L)}
 \cong
 (\WW\sporb G)_{(L)}\\
 &= 
 T^*(Q_{(H)}/G)\times_{Q_{(H)}/G}((\bsc_{q\in Q_{(H)}}\orb\cap\ann\gu_{q})_{L})/G
\end{align*}
are computed by the same formulas.
\end{enumerate}
\end{theorem}

\begin{proof}
Property (1) follows directly from 
the observation that 
\begin{align*}
 (T^*Q_{(H)})\sporb G
 &\cong
 T^*(Q_H\times_{W}G/H)\sporb G\\
 &\cong 
 (T^*Q_H\times(T^*G)\spr{0}H)\spr{0}W\sporb G
 =
 (T^*Q_H\times\orb\spr{0}H)\spr{0}W
\end{align*}
where we use commuting reduction to get the second \iso which thus is
strata and Poisson structure preserving. 
The first \iso of this equation is stratified and Poisson structure
preserving since it is constructed from the lifting of the
$G$-equivariant diffeomorphism $Q_{(H)}\cong Q_H\times_WG/H$.
Hereby, as above, $H$ and $W$ act on $G$ by inversion of right
multiplication, and these actions are cotangent lifted to $T^*G$. 

(2)
Indeed, $W$ acts freely and by Hamiltonian transformations on
$T^*Q_H\times\orb\spr{0}H$. The \momap of this action is 
\[
 J_W: T^*Q_H\times\orb\spr{0}H\longto\wo^*,
 (q,p,[\lam]_H)\longmapsto\mu(q,p)-\lam|\wo
\]
which is well-defined and computed by the same reduction-in-stages
argument as in Lemma \ref{lem:J_W}.
Let $\mathcal{A}^*: Q_H\times\wo^*\to\bsc_{q\in Q_H}T_q^*(W.q) = \ver^*$,
$(q,\lam)\mapsto\mathcal{A}^*_q(\lam)$ be the dual of $\mathcal{A}$.
By construction $\mathcal{A}_q^*$ is an inverse to
$\mu_q|\ver_q^*$. Therefore,
\[
 J_W^{-1}(0)
 \cong_{\mathcal{A}}
 \set{(q,p_0+\mathcal{A}_q^*(\lam|\wo),[\lam]_H):
      p_0\in\hor_q^*,[\lam]_H\in\orb\spr{0}H}
 \cong
 \hor^*\times\orb\spr{0}H.
\]
which is an \iso of stratified spaces since $\mu_q$ clearly is
$H$-equivariant. Thus 
\[
 (T^*Q_H\times\orb\spr{0}H)\spr{0}W
 =
 J_W^{-1}(0)/W
 \cong_{\mathcal{A}}
 \hor^*\times_W\orb\spr{0}H\longto 
 \hor^*/W=T^*B
\]
is the fiber bundle over $T^*B$ associated to the Hamiltonian
$W$-action on $\orb\spr{0}H$. 

(3)
Writing $Q_{(H)}$ as an associated bundle we see, as above, that 
$T^*Q_{(H)}\cong 
 (T^*Q_H\times T^*G\spr{0}H)\spr{0}W$ as smooth symplectic
manifolds. Since $H$ is normal in $N(H)$ by tautology  
the Regular Reduction in Stages Theorem
(\cite{MMOPR03,OR04}) implies
that 
$T^*Q_{(H)}\cong(T^*Q_H\times T^*G)\spr{0}N$ where $N:=N(H)$. 
Therefore, we can describe the symplectic stratum $(T^*Q_{(H)}\sporb
G)_{(L)} := (\mu^{-1}(\orb)\cap(T^*Q_{(H)})_{(L)})/G$ as follows.
\begin{align*}
 (T^*Q_{(H)}\sporb G)_{(L)}  
 &\cong
 ((T^*Q_H\times T^*G)\spr{0}N\sporb G)_{(L_0)}\\
 &\cong\tag{2}\label{eq:2}
 ((T^*Q_H\times\orb)\spr{0}N)_{(L_0)^N}\\
 &=
 (J_N^{-1}(0)\cap(T^*Q_H\times\orb)_{(L)^N})/N\\
 &= 
 (J_N^{-1}(0)\cap(T^*Q_H\times\orb_{(L)^N}))/N\\
 &=
 \set{(q,p,\lam):\lam|\no=\mu(q,p)\in\wo^*\subset\ann\ho}/N\\
 &=
 \set{(q,p,\lam):\lam\in\orb\cap\ann\ho, \lam|\wo=\mu(q,p)}/H/W\\
 &=
 \set{(q,p,[\lam]_H)\in J_W^{-1}(0)\subset T^*Q_H\times(\orb\spr{0}H)_{(L_0)^N}}/W\\
 &=\tag{8}\label{eq:8}
 (T^*Q_H\times W.(\orb\spr{0}H)_{(L_0)^H})\spr{0}W
\end{align*}
where $J_N$ is the $N$-\momap on $T^*Q_H\times\orb$ given by
$J_N(q,p,\lam) = \mu(q,p)-\lam|\no$, $\no$ is the Lie algebra of $N$,
and the non-obvious identifications are verified as follows.
For identification (\ref{eq:2}) notice that $N\times G$ acts in a
Hamiltonian and proper fashion on $T^Q_H\times T^*G$. A typical
isotropy group of this action is of the form
\[
(N\times G)_{(q,p,k,\lam)}
 =
 \set{(h,khk^{-1}): h\in H_{\lam}} =: L'.
\]
It is straightforward to check that $((0,\lam_0),L')$ (where
$\lam_0\in\orb\cap\ann\ho$) satisfies the Hamiltonian Stages
Hypothesis for non-free actions as formulated in 
\cite[Section~10.4]{MMOPR03} or in \cite[Section~9.5]{OR04}. Thus we
can apply the corresponding reduction theorem as given in these
references. The result is identification (\ref{eq:2}). As stated in
Remark~\ref{rem:conn} the (induced) momentum maps appearing in these
computations are equivariant whence we are not concerned with the
connectedness hypothesis that is made in the general formulation of
the Reduction in Stages Theorem. 
The last identification (\ref{eq:8}) is a consequence of 
Lemma~\ref{lem:W-sweep}.

Item (4)
is an application of the result in Sternberg \cite{S77}.  
See also \cite{W78,GLS96}.

\begin{comment}
***********************************************************************
VERRECHNET UND UNNOETIG!!!!!!!!
In order to see property (4) we suppress the obvious pull-backs by
$\tau$, $\pr_1$, and $\pr_2$, and let $\psi=\psi(\mathcal{A})$ 
denote the \iso 
\begin{align*}
 \hor^*\times(\orb\spr{0}H)_{(L_0)^H}&\longto J_W^{-1}(0)\subset T^*Q_H\times(\orb\spr{0}H)_{(L_0)^H},\\
 (q;b,p_0,[\lam])&\longmapsto(q,(T_q\pi_W)^*p_0+\mathcal{A}_q^*(\pr_0[\lam]),[\lam])
\end{align*}
where $\pr_W: Q_H\toto B$ is the orbit projection.
(See Proposition \ref{prop:cot-red}.) 
Let us denote the canonical
symplectic structure on $T^*Q_H\times(\orb\spr{0}H)_{(L_0)^H}$ by
$\Om^{Q_H}+\lorb{\Om} = -d\theta^{Q_H}+\lorb{\Om}$, and for notational
simplicity we let its restriction to $J_W^{-1}(0)$ go by the same
name. 
Since $\psi^*\theta^{Q_H}$ at a point $(q;b,p_0,[\lam];q',b',p_0',[\lam]')\in
T(\hor^*\times(\orb\spr{0}H)_{(L_0)^H})$ evaluates to
\begin{multline*}
 \theta^{Q_H}(q,(T_q\pi_W)^*p_0+\mathcal{A}_q^*(\pr_0[\lam]))(q',T^v_{(q;b,p_0,[\lam])}.(q',b',p_0',[\lam]'))
 =\\
 = \vv<p_0,T_q\pi_W.q'> + \vv<\pr_0[\lam],\mathcal{A}_q(q')>
\end{multline*}
we see that $\psi^*\theta^{Q_H} =
\eta^*\theta^{B}+\vv<\pr_0,\mathcal{A}>$, whence the assertion
follows. 
(Hereby,
$T^v_{(q;b,p_0,[\lam])}.(q',b',p_0',[\lam]')$ denotes the -irrelevant-
vertical part of the tangent vector 
$T_{(q;b,p_0,[\lam])}.(q',b',p_0',[\lam]')$.)
**************************************************************************
\end{comment}

Assertions (5) follows from
Proposition \ref{prop:curv}, i.e., from the fact that $\curv^A$ can be
computed from $\curv^{\mathcal{A}}$. (For the formula determining the
symplectic structure on $\WW\sporb G$ see \cite[Theorem 5.5]{H04}.)
\end{proof}

\begin{rem}
Notice that this theorem allows to compute the reduced cotangent
bundle $T^*Q_{(H)}\sporb G$ without explicitly knowing the \momap 
$\mu: T^*Q_{(H)}\to\gu^*$. 
\end{rem}

\begin{rem}
Looking again at the proof of Item \ref{item:iso} in the above theorem
one could also directly apply the singular version of the Reduction in
Stages Theorem as follows. Namely, $N$ acts in a Hamiltonian fashion
on $T^*Q_H\times\orb$, a typical isotropy type of this action is
$(L_0)^N$, and $H$ is normal in $N$. One can check again that the pair
$(0,L_0)$ satisfies the Hamiltonian Stages Hypothesis. Thus the
Reduction in Stages Theorem is applicable and using it we get a
symplectomorphism
\[
 (T^*Q_{(H)})\sporb G)_{(L)}
 \cong
 (T^*Q_H\times(\orb\spr{0}H)_{(L_0)^H})\spr{0}\by{N_{N(H)}(L_0)}{N_H(L_0)}
\]
which is an equivalent description of the strata of the reduced
space. By Lemma~\ref{lem:strat-pres} these descriptions coincide when
$(L_0)^N=(L_0)^H$.
\end{rem}

\begin{rem}[Geometric quantization]\label{rem:quant}
According to the above theorem a prescribed connection $\mathcal{A}$ 
yields a symplectic fibration of 
$T^*Q_{(H)}$
over $T^*B$ with fiber
$\orb\spr{0}H$. Using the quantization in stages procedure outlined in
\cite[Section~4.1]{GLS96} one could thus try to quantize
$T^*Q_{(H)}\sporb G$ via quantization of base $T^*B$ and fiber
$\orb\spr{0}H$. It might be interesting to investigate whether this has
any useful consequences for the quantization of the (singular)
symplectic quotient $T^*Q\sporb G$.
\end{rem}

We end this section by exploring some particular cases of Theorem
\ref{thm:red}.

\subsubsection{A lot of symmetry}
Suppose $L_0=H$ such that $(\orb\spr{0}H)_{(H)^H} = \orb_H\cap\ann\ho
= \orb\cap\wo^*$ which obviously is invariant under $W$. Then,
the smooth minimal coupling space
\[ 
 ((T^*Q_{(H)})\sporb G)_{(H)}
 \cong
 (T^*Q_H\times(\orb\cap\wo^*))\spr{0}W 
 \cong_{\mathcal{A}}
 (Q_H\times_B T^*B)\times_W\orb\cap\wo^*
\]
is a Hamiltonian fiber bundle over $T^*B$.

\subsubsection{A lot of regularity}
Assume $\orb$ is compact.
By Theorem \ref{thm:open-dense} there is a unique maximal conjugacy
class of the $H$-action on $\orb$. Suppose $(L_0)^H$ is this
class. Then $(\orb\spr{0}H)_{(L_0)^H}$ is -as is noted in Lemma
\ref{lem:W-sweep}- invariant under the $W$-action. 
Therefore,
the smooth minimal coupling space
\begin{align*}
 ((T^*Q_{(H)})\sporb G)_{(L_0)}
 &\cong
 (T^*Q_H\times(\orb\spr{0}H)_{(L_0)^H})\spr{0}W\\ 
 &\cong_{\mathcal{A}}
 (Q_H\times_B T^*B)\times_W(\orb\spr{0}H)_{(L_0)^H}
\end{align*}
is a Hamiltonian fiber bundle over $T^*B$.

\subsubsection{Hamiltonian reduction}
Assume $G$ acts on $Q$ by isometries with respect to a Riemannian
structure $\vv<.,.>$.
Let $\mathcal{H}$ be the free Hamiltonian on $T^*Q$ given by the
metric on $Q$, that is $\mathcal{H}(q,p)=\by{1}{2}\vv<p,p>$. 
To simplify the notation we use the same symbol $\vv<.,.>$ for the
metric as well as for the cometric.
Then
Theorem \ref{thm:red} gives a way of computing the reduced Hamiltonian
$\mathcal{H}_0$ on 
the reduced phase space
$(T^*Q_{(H)})\sporb G \cong
(T^*Q_H\times\orb\spr{0}H)\spr{0}W$. Namely, for $q\in Q_H$ let
$\ine^G_q$ be the non-degenerate pairing on $\ho^{\bot}$ defined by
$\ho^{\bot}\times\ho^{\bot}\to\R$,
$(X,Y)\mapsto\vv<\zeta_X(q),\zeta_Y(q)>$. Note that this pairing can
be extended $G$-equivariantly to define a tensor $\ine^G:
Q_{(H)}\to\bsc_{q\in Q_{(H)}}(\gu/\gu_q)^*\otimes(\gu/\gu_q)^*$. This
tensor is used in \cite{H04} to define the generalized mechanical
curvature $A$ on the bundle $Q_{(H)}\toto Q_{(H)}/G$. Notice further,
that the inertia tensor $\ine=\ine^W: Q_H\to\wo^*\otimes\wo^*$ 
defined in Section \ref{sec:gen}
is just
the restriction of $\ine^G$ to $\bsc_{q\in Q_H}\wo^*\otimes\wo^*$. 
Conversely, it is not true that $\ine^G$ can be computed from
$\ine^W$. However, by virtue of Proposition \ref{prop:curv} we can
obtain the generalized mechanical curvature $\curv^A$ from
the simpler form $\curv^{\mathcal{A}}$. 

Regarding the computation of $\mathcal{H}_0$ let
$\lam\in\ann\ho=(\gu/\ho)^*$, and let $X_q(\lam)$ denote the vector in
$\ho^{\bot}$ determined by the pairing $\ine^G_q$. Then the reduced
Hamiltonian is given by
\begin{align*}
 &(T^*Q_H\times\orb\spr{0}H)\spr{0}W
 \cong_{\mathcal{A}}
 (Q_H\times_B T^*B)\times_W\orb\spr{0}H
 \longto\R\\
 &[(q;q_0,p_0;[\lam]_H)]_W
 \longmapsto
 \by{1}{2}\vv<p_0,p_0>_{q_0}+\by{1}{2}\ine_q^G(X_q(\lam),X_q(\lam)).
\end{align*}
This follows immediately from the identity
$\vv<{A}_q^*(\lam),{A}_q^*(\lam)> = 
\ine_q^G(X_q(\lam),X_q(\lam))$ where ${A}_q^*$ is the point
wise dual to ${A}_q$ determined by the metric and the inertia
pairing on $\ho^{\bot}$. 

As the reduced phase space is a stratified symplectic space we get 
via restriction of $\mathcal{H}_0$
a Hamiltonian system in the usual sense on each stratum 
$((T^*Q_{(H)})\sporb G)_{(L_0)} \cong
(T^*Q_H\times W.(\orb\spr{0}H)_{(L_0)^H})\spr{0}W$
where $L_0\subset H$.

\begin{rem}[Wong's equations]\label{rem:wong}
Suppose $\beta$ is a $G$-biinvariant metric on $\gu$, and let $g$ be a
$G$-invariant Riemannian metric on $Q=Q_{(H)}$ such that $\ine_q^G=\beta$,
independently of $q\in Q$. The corresponding free Hamiltonian system
on $T^*Q$ is the Kaluza-Klein system of Kerner~\cite{Ker68}. In the
case that the $G$-action on $Q$ is free, Montgomery~\cite{Mon84} has
shown that Hamiltonian reduction of the Kaluza-Klein system yields
Wong's equations (Wong~\cite{Won70}). It is further observed in
\cite{Mon84} that the reduced Kaluza-Klein system is equivalent to
Sternbergs minimal coupling Hamiltonian on  
$T^*Q\sporb G \cong (T^*Q_H\times\orb\spr{0}H)\spr{0}W
 \cong_{\mathcal{A}}
 T^*B\times_B(Q\times_W\orb\spr{0}H)
 \to T^*B$,
whence it is also equivalent to Weinstein's description of the reduced
Hamiltonian system. (See \cite{S77,W78} or
Subsection~\ref{sub:min}). These results carry over to the singular
situation as well. 
This is roughly seen as follows.
Using the Slice Theorem to get a local description
of $T^*Q\sporb G$ as in \cite[Theorem~4.4]{H04} one can mimic the
computations of \cite{Mon84} to obtain Wong's equations from the
Kaluza-Klein Hamiltonian system. 
By Theorem~\ref{thm:red}(5) the resulting Hamiltonian system is
equivalent to the system obtained via Hamiltonian reduction of
$(T^*Q_H\times\orb\spr{0}H,
 \pr_1^*\Om^{Q_H}+\pr_2^*\gamma,
 \by{1}{2}g+\by{1}{2}\beta)$
where
$(\pr_1,\pr_2): T^*Q\times\orb\spr{0}H\to T^*Q\times\orb\spr{0}H$
denote the Cartesian projections 
and $\gamma$ is the stratified form on $\orb\spr{0}H$ which restricts
to the induced symplectic form on each of the smooth symplectic strata
of $\orb\spr{0}H$.
The resulting Wong's equations for a curve $(c(t),p(t),\lam(t))$ on a
stratum of $T^*B\times_B(Q_H\times_W\orb\spr{0}H) \cong{\mathcal{A}}
T^*Q_{(H)}\sporb G$ can be stated as 
\[
 \nabla_{c'}c'
 = 
 -\hat{g}^{-1}(\vv<\lam,i_{c'}\curv_0^{\mathcal{A}}>)
 \text{ and }
 D_{c'}\lam = 0
\]
where $D_{c'}$ denotes covariant differentiation 
along $c$
of sections of the 
bundle $Q_H\times_W\orb\spr{0}H$ with respect to the
connection associated to $\mathcal{A}$. This is well defined on each
of the smooth strata of this bundle. (The form of these equations is
the same as that in Montgomery~\cite{Mon02} where also a general
discussion of Wong's equations can be found.)
Moreover by the above, the reduced
Kaluza-Klein Hamiltonian is of the form 
\[
 (Q_H\times_B T^*B)\times_W\orb\spr{0}H
 \longto\R,
 [(q;q_0,p_0;[\lam]_H)]_W
 \longmapsto
 \by{1}{2}\vv<p_0,p_0>_{q_0}+\by{1}{2}\beta(\lam,\lam)
\]
where we use the same symbol $\beta$ for the dual metric on
$\gu^*$. As in the regular case this Hamiltonian differs from
Sternberg's minimal coupling Hamiltonian only by a Casimir
function. Thus the resulting equations of motion coincide. 
\end{rem}

\subsubsection{Another set of Hamiltonian equations}
We continue to assume that $G$ acts on $(Q,\vv<.,.>)$ by
isometries. Thus $W$ acts on $(Q_H,\vv<.,.>)$ by isometries as well,
and reduction at $0$ with respect to $W$ of the free
Hamiltonian system $(T^*Q_H,\Om^{Q_H},\by{1}{2}\vv<.,.>)$ 
yields the free Hamiltonian system on
$T^*B$.  
Let $(\pr_1,\pr_2): T^*Q\times\orb\spr{0}H\to T^*Q\times\orb\spr{0}H$
denote the Cartesian projections, and reduce the $W$-invariant system
$(T^*Q_H\times\orb\spr{0}H,\pr_1^*\Om^Q+\pr_2^*\gamma,\pr_1^*\by{1}{2}\vv<.,.>)$ 
at $0\in\wo^*$
with respect to the diagonal $W$-action. The Hamiltonian of
this reduced system is given by 
\begin{align*}
 &(T^*Q_H\times\orb\spr{0}H)\spr{0}W
 \cong_{\mathcal{A}}
 (Q_H\times_B T^*B)\times_W\orb\spr{0}H
 \longto\R\\
 &[(q;q_0,p_0;[\lam]_H)]_W
 \longmapsto
 \by{1}{2}\vv<p_0,p_0>_{q_0}+\by{1}{2}\ine_q^W(X_q(\lam|\wo),X_q(\lam|\wo)).
\end{align*}
Comparing this expression to the above it is clear that this describes
a Hamiltonian system which is in general different from that obtained
by Hamiltonian reduction at $\orb$ of $(T^*Q_{(H)},\Om,\by{1}{2}\vv<.,.>)$.
However, 
considering the restriction of the reduction of the free system on
$T^*Q_{(H)}$ to the stratum 
$((T^*Q_{(H)})\sporb G)_{(H)}
 \cong
 (T^*Q_H\times(\orb\cap\wo^*))\spr{0}W 
 \cong_{\mathcal{A}}
 (Q_H\times_B T^*B)\times_W\orb\cap\wo^*$
shows that, on this stratum, these a priori different systems
coincide. Also, if the $G$-action on $Q_{(H)}$ is free these systems
coincide by virtue of the shifting trick.

\subsection{Appendix: Minimal coupling and symplectic fibrations}\label{sub:min}
The purpose of this appendix is to shortly say what we mean by minimal
coupling. A detailed exposition of the subject can be found in
Guillemin, Lerman, and Sternberg \cite{GLS96} which is also the
reference for the subsequent. 
Let $G\hookto Q\toto B$
be a principal fiber bundle equipped with a principal bundle
connection form $\mathcal{A}\in\Om^1(Q;\gu)$, and suppose $(F,\Om^F)$
is a right Hamiltonian $G$-space with \momap $J_F: F\to\gu^*$. In \cite{S77}
Sternberg has shown how to construct from these data a symplectic form
on the so-called Sternberg space
\[
 (Q\times_B T^*B)\times_G F.
\] 
Weinstein \cite{W78} noticed that this result can be obtained in a
more symplectic way through the following universal procedure. Namely,
cotangent lift the action by $G$ on $Q$ 
(which we assume to be a proper left action)
to the cotangent bundle $T^*Q$. 
Thus we can 
consider the diagonal action on $T^*Q\times F$
where the action on $F$ is inverted.
This action is
Hamiltonian with \momap $J := \mu-J_F$ -- where $\mu$ is the cotangent
bundle \momap. Thus we can do symplectic reduction to get a new symplectic
manifold $J^{-1}(0)/G = (T^*Q\times F)\spr{0}G$. Choosing now a
connection $\mathcal{A}$ on $Q\toto B$ yields an explicit
symplectomorphism
with Sternberg's $\mathcal{A}$-dependent space, i.e.,
\[
 (T^*Q\times F)\spr{0}G
 \cong_{\mathcal{A}}
 T^*B\times_B(Q\times_G F)
 \cong
 (Q\times_B T^*B)\times_G F.
\]
This point of view is generally called the Weinstein space picture. 
Either of the isomorphic ways of
constructing a symplectic fiber bundle (definition below) 
$\chi=\chi(\mathcal{A}): (T^*Q\times F)\spr{0}G\to T^*B$
out of the data $G\hookto
Q\toto B$, $\mathcal{A}$, and $(F,\Om^F)$ is referred to as
\caps{minimal coupling}. 

The idea is that one may thus start from a
Hamiltonian system $(T^*B,\Om^B,\mathcal{H})$ and obtain a new
system on
$(T^*Q\times F)\spr{0}G$, 
together with its induced symplectic structure,
with respect to the Hamiltonian
$\eta^*\mathcal{H}$. 
This way of producing a new
system is thought of a adjoining some sort of internal variables (like
spin) to the original system on $B$, and the connection $\mathcal{A}$
is interpreted as the potential of a Yang-Mills field which affects the
system on $B$. 

A \caps{symplectic fiber bundle} is a fiber bundle $F\hookto X\xrightarrow{\pi}M$
with fiber $(F,\Om^F)$ a symplectic manifold 
such that the transition functions take values in the group of
symplectomorphisms of $F$. Given such data one may ask whether there
is a symplectic form on $X$ such that its restriction to a fiber is
the prescribed symplectic form on the fiber? In general there is no
such globally defined form on $X$. However, one can give a partial
affirmative answer to this question via the coupling form: A two-form
$\om$ on $X$ is called \caps{fiber compatible} if its restriction to a
fiber is the prescribed symplectic form on that fiber. Let $\ver :=
\ker T\pi$ and suppose $\Gamma$ is a connection on $\pi: X\to M$ such
that we also have a horizontal subbundle $\hor=\hor(\Gamma)$. 
The connection
$\Gamma$ is called \caps{symplectic} if the associated parallel
transport is fiber-wise symplectomorphic. Given such a connection one
can define a fiber compatible form $\om(\Gamma)$ by declaring
$i_v\om(\Gamma)=0$ for all horizontal fields $v$ and letting
$\om(\Gamma)$ restrict to the prescribed symplectic form on each
fiber. Conversely, given a fiber compatible form $\om$ on $X$ we can
define
\[
 \hor(\om)
 :=
 \set{v\in TX: \om(v,w)=0 \text{ for all } w\in\ver}
 :=
 \ver^{\om}.
\]
If $\hor(\om)=\hor(\Gamma)$ then $\om$ is said to be
$\Gamma$-\caps{compatible}. 
This definition does not depend on $\Gamma$ being symplectic.
Clearly, $\om(\Gamma)$ is
$\Gamma$-compatible. 

There exist strong results 
(\cite{GLS96,GLSW83})
concerning two-forms on the total space of
a symplectic fiber bundle $F\hookto X\xrightarrow{\pi}M$.
In particular it is true that every symplectic fiber bundle has a
symplectic connection. Since we
will not make use of these results we refrain from stating them
explicitly and just refer to \cite{GLS96,GLSW83} which contain a
detailed discussion. 
As a matter of fact, \cite[Theorem~1.4.1]{GLS96} provides a way of
assuring existence and uniqueness after suitable normalization of
$\Gamma$-compatible forms for the case that $F$ is compact, connected,
and simply connected. This uniquely characterized form is called
\caps{minimal coupling form} of the symplectic fibration.

\begin{comment}
We now reproduce two remarkable results
concerning two-forms on the total space of the
symplectic fiber bundle $F\hookto X\xrightarrow{\pi}M$. Both 
can be found in \cite{GLS96}.
\begin{theorem}[\cite{GLSW83}]
Let $\Gamma$ be a connection on $F\hookto X\xrightarrow{\pi}M$, and
let $\om$ be a $\Gamma$-compatible two-form on $X$. Then $\Gamma$ is
symplectic if and only if $i_{v\wedge w}d\om = 0$ for all vertical
vector fields $v,w$. 
\end{theorem} 
As a corollary of this theorem one easily sees that every symplectic
fiber bundle has a symplectic connection. 
Indeed, this is proved using a partition of unity argument over $M$.
\begin{theorem}[\cite{GLS96}]
Suppose $F$ is compact, connected, and simply connected, and assume
that $\Gamma$ is a symplectic connection $\pi: X\to M$. Then there
exists a unique closed $\Gamma$-compatible two-form $\om_{\Gamma}$
with the property that $\pi_*\om_{\Gamma}^{d+1} = 0$ where $\pi_*$ is
the fiber integration map.  
\end{theorem}
The form $\om_{\Gamma}$ is called the \caps{minimal coupling form}
associated to the symplectic connection $\Gamma$. 
\end{comment}

Let us return to the minimal coupling construction above. 
The minimal coupling form of $\pr_1: T^*Q\times F\to T^*Q$ is simply
$\pr_2^*\Om^F$. 
Further,
let $\mathcal{A}$ continue to denote the principal
bundle connection form on $G\hookto Q\toto B$. 
This gives rise to a
symplectic connection $\Gamma=\Gamma(\mathcal{A})$ on 
$F\hookto X=(T^*Q\times F)\spr{0}G\to T^*B = M$. 
Concerning the associated
minimal coupling form $\om_{\Gamma}$ let $\tau: T^*Q\to Q$ be the
footpoint projection, and consider
$\pr_2^*\Om^F-\vv<\pr_2,(\tau\circ\pr_1)^*\curv^{\mathcal{A}}>$. This
form restricts to an horizontal and $G$-invariant object on
$J^{-1}(0)$ and drops to $\om_{\Gamma}$ via the orbit projection
$J^{-1}(0)\toto J^{-1}(0)/G = X = (T^*Q\times F)\spr{0}G$.
Again, for details we refer to \cite{GLS96}.

A \caps{Hamiltonian fiber bundle} is a symplectic fiber bundle whose
transition functions take values in the group of 
Hamiltonian transformations of the fiber. 

Therefore, our interpretation of Theorem \ref{thm:red} as exhibiting
$(T^*Q_{(H)})\sporb G$ as a minimal coupling space is justified.

\section{Remarks on the stratification of cotangent bundles}\label{sec:strata}
Let $(H)$ denote an isotropy type of the $G$-action on $Q$, and $(L)$
be an isotropy type of the cotangent lifted $G$-action.
One of the obvious problems with singular reduction of $T^*Q$ with
respect to $G$ is that the foot point projection $\tau: T^*Q\to Q$ is
not stratified, i.e., the preimage of a stratum under $\tau$ is not
equal to a union of strata of $T^*Q$. This poses a problem for the
Hamiltonian dynamics. Indeed, the
Hamiltonian flow of a $G$-invariant function is easily seen to
preserve strata $(T^*Q)_{(L)}$, however, it generally neither
preserves $(T^*Q)_{(L)}|Q_{(H)}$ nor $(T^*Q)|Q_{(H)}$. Forcing $\tau$
to be stratified by further decomposing strata $(T^*Q)_{(L)}$ into
pieces of the form $(T^*Q)_{(L)}|Q_{(H)}$ obviously yields a finer
stratification of $T^*Q$ which thus needs to be studied. This finer
decomposition of $T^*Q$ was called \caps{secondary stratification} in
Perlmutter et al.\ \cite{PRS03}, and we shall adopt this terminology.

Let $\ann Q_{(H)}\to Q_{(H)}$ denote the subbundle of $(T^*Q)|Q_{(H)}$
consisting of those covectors which vanish upon insertion of a vector
tangent to $Q_{(H)}$. Clearly, we have
\[
 (T^*Q)_{(L)}|Q_{(H)}
 =
 (T^*Q_{(H)}\times_{Q_{(H)}}\ann Q_{(H)})_{(L)},
\]
and note that the \momap $\mu: T^*Q\to\gu^*$ vanishes on 
$\ann Q_{(H)}$. 
Therefore,
for an orbit $\orb$ in the image of $\mu$ we have that 
\[
 \mu^{-1}(\orb)|Q_{(H)}
 = 
 \mu_{(H)}^{-1}(\orb)\times_{Q_{(H)}}\ann Q_{(H)}
\]
where $\mu_{(H)}$ denotes the \momap of the cotangent lifted
$G$-action on $T^*Q_{(H)}$. 
Thus one should study the fibration
\[
 (\ann_q Q_{(H)})/L\hookto(\mu_{(H)}^{-1}(\orb)\times_{Q_{(H)}}\ann Q_{(H)})_{(L)}/G
 \longto(\mu_{(H)}^{-1}(\orb))_{(L)}/G.
\]
If $Q_{(H)}=Q_{\textup{reg}}$ is the regular stratum which is open dense in
$Q$ then $\ann Q_{(H)}$ is trivial. Thus in this case
Theorem~\ref{thm:red} gives a full answer to the reduction problem,
and this is the generic case.  

Suppose now that $Q_{(H)}=Q_{\textup{reg}}$ 
is the regular stratum of $Q$. Similarly as in Proposition
\ref{prop:cot-red} we see that 
\begin{align*}
 (T^*Q)_{(L)}|Q_{(H)}
 &=
 (T^*Q_{(H)})_{(L)}
 \cong
 (T^*(Q_H\times_{W}G/H))_{(L)}\\
 &=
 (T^*Q_H\times G\times_H\ann\ho)_{(L)}\spr{0}W\\
 &\cong
 (Q_H\times_{B}T^*B)\times_W(G\times_H\ann(\ho+\wo)\times\wo^*)_{(L)}
\end{align*}
where $B=Q_H/W$. The problem in this case thus reduces to understanding the
$G$-action on $G\times_H\ann(\ho+\wo)$ which is given by $g.[(k,\lam)]
= [(gk,\lam)]$. 
Now,
$[(k,\lam)]\in(G\times_H\ann(\ho+\wo))_{(L)}$ if and only if
$k^{-1}G_{[(k,\lam)]}k = G_{[(e,\lam)]} =: L_0 \sim L$, and the latter
is the case if and only if $\lam\in(\ann\ho)_{L_0}$ with respect to
the $H$-action on $\ann\ho$. In particular, $L_0\subeq H$.

\section{Some well hidden symmetries}\label{sec:hidsym}
Let the matrix group $G:=SO(5)$ act on $Q:=S^9\subset\R^5\times\R^5=\R^{10}$ through the
diagonal action. With respect to this action $Q$ decomposes into two
orbit type strata corresponding to $H_0:=SO(3)\hookto G$ and
$H_1:=SO(4)\hookto G$. Both embeddings are the standard embeddings into the
lower right corner of the matrices in $G$. Using minus one half of 
the trace form on
$\gu=\so(5)$ we shall tacitly identify $\gu$ and $\gu^*$.

\subsection{The orbit space $S^9/SO(5)$}
We write elements $q\in Q\subset\R^5\times\R^5$ as
$q=(q_1^j,q_2^j)_{j=1}^5$. 
The subset of regular elements $Q_{(H_0)}$ is the set of $(q_1,q_2)\in
Q$ such that $q_1$ and $q_2$ are linearly independent. In view of the
isomorphism $Q_{(H_0)}/G\cong Q_{H_0}/W_0$ where $W_0:=N(H_0)/H_0$
(and Proposition \ref{prop:cot-red}) we need to identify
$Q_{H_0}$ and $W_0$. With the embedding  $H_0=SO(3)\hookto G$ into the lower
right corner we have
\begin{multline*}
 Q_{H_0}
 =
 \{((a,b,0,0,0)^t,(\alpha,\beta,0,0,0)^t):\\
      (a,b)^t \text{ and } (\alpha,\beta)^t 
       \text{ are linearly independent and }
       a^2+b^2+\alpha^2+\beta^2 = 1\}
\end{multline*}
where $(\phantom{x})^t$ denotes transpose. We embed $Q_{H_0}\hookto
S^3$ in the obvious way.
Moreover,  
\[
 W_0
 =
 S^1\times\{
  \pm\left(
   \begin{matrix}
    1&0\\
    0&1
   \end{matrix}
  \right),
  \pm\left(
   \begin{matrix} 
    1&0\\
    0&-1
   \end{matrix}
  \right)
  \} 
 =:
 S^1\times\Delta
\]
acts on $Q_{H_0}\subset S^3$ in the standard fashion.
Thus the orbit
projection $\eta_0: Q_{H_0}\toto Q_{H_0}/S^1$ is just a
restriction of the Hopf map to an open dense subset 
$Q_{H_0}\subset S^3$. 
In order to explicitly describe $B_0=Q_{H_0}/W_0$ 
it is convenient to write the
Hopf map $\eta$ from $S^3$ onto the sphere $S^2(\by{1}{2})$ of radius 
$\by{1}{2}$ as 
\begin{align*}
 \eta: \R^4=\C^2\supset S^3&\longto S^2(\by{1}{2})\\
 (q_1,q_2)=(a+ib,\alpha+i\beta)
  &\longmapsto
 (\by{1}{2}(|q_1|^2-|q_2|^2),\Re{(q_1\ov{q_2})},\Im{(q_1\ov{q_2})})^t.
\end{align*}
Now it is straightforward to see that $\eta_0=\eta|Q_{H_0}$ takes
$Q_{H_0}$ onto the open subset
$\set{(x,y,z)^t\in S^2(\by{1}{2}): z\neq0}$, and the $\Delta$-action
factors through $\eta_0$ to an action by $\set{\pm1}$ in the
$z$-direction. Therefore, 
$B_0 = Q_{H_0}/W_0 = \set{(x,y,z)^t\in S^2(\by{1}{2}): z>0}$, and this
is the regular stratum of the orbit space $B:=Q/G$. 

Similarly, to get the other stratum of $B$ notice that 
\[
 Q_{H_1}
 =
 \{((a,0,0,0,0)^t,(\alpha,0,0,0,0)^t):
       a^2+\alpha^2 = 1\}
\]
and $W_1:=N(H_1)/H_1$ is trivial. Thus $Q_{H_1}/W_1\cong S^1$, and via
$\eta$ we can send this $S^1$ diffeomorphically onto the equator in
the $x$-$y$-plane in $S^2(\by{1}{2})$. Therefore,
$B = \set{(x,y,z)^t\in S^2(\by{1}{2}): z\ge0} = 
B_0\sqcup B_1$ is stratified into northern hemisphere plus equator.

\subsection{The reduced phase space $T^*S^9/SO(5)$}
Now we can invoke Proposition \ref{prop:cot-red} to compute
$(T^*Q_{(H_0)})/G \cong
(Q_{H_0}\times_{B_0}T^*B_0)\times_{W_0}(\wo_0^*\times\ann(\ho_0+\wo_0)/H_0)$. 
Indeed, $W_0$ acts trivially on $\wo_0 = \R$ whence the
reduced phase space is an associated bundle (in a singular sense --
see \cite{hr05}) of the type
\begin{multline*}
 \R\times\ann(\ho_0+\wo_0)/H_0\cong
 \R\times\R^3\times_{H_0}\R^3\hookto\\
 \hookto(T^*Q_{(H_0)})/G
 \cong
 (Q_{H_0}\times_{B_0}T^*B_0)\times_{W_0}(\R\times\R^3\times_{H_0}\R^3)
 \longto T^*B_0
\end{multline*}
where we use the negative of the trace form on $\gu$ to identify
$\ann(\ho_0+\wo_0)\subset\gu^*$ and
$\ho_0^{\bot}\cap\wo_0^{\bot}\subset\gu$ 
as well as the
$H_0$-equivariant linear \iso
\[
 \ho_0^{\bot}\cap\wo_0^{\bot}\longto\R^3\times\R^3,
 \text{ }
 (x_{ij})_{ij}\longmapsto((x_{k1})_{k=3}^{5},(x_{k2})_{k=3}^{5})
  = (v,w)
\]
where $H_0=SO(3)$ acts on $\R^3\times\R^3$ in the standard diagonal
way.

According to the remarks in Section \ref{sec:strata} we also want to
consider the regular stratum $(T^*Q)_{reg}=(T^*Q)_{(L_0)}$ as well as
its secondary strata $(T^*Q)_{(L_0)}|Q_{(H_0)}$ and
$(T^*Q)_{(L_0)}|Q_{(H_1)}$. It is easy to see that $L_0=\set{1}$ and that
$(T^*Q)_{(L_0)}|Q_{(H_1)}$ is empty. Furthermore,
$(G\times_H\ann(\ho_0+\wo_0))_{(L_0)} \cong
(G\times_H(\R^3\times\R^3))_{\{1\}}$ consists of those elements
$[(k,v,w)]$ for which $v$ and $w$ are linearly independent, that is
$(v,w)\in(\R^3\times\R^3)_{\{1\}}$. Therefore, the regular part of
the reduced phase space $(T^*Q)/G$ is
\begin{multline*}
 \R\times\R_{>0}\times\R\times\R_{>0}\hookto\\
 \hookto(T^*Q_{(H_0)})_{\{1\}}/G
 \cong
 (Q_{H_0}\times_{B_0}T^*B_0)\times_{W_0}(\R\times\R^3\times_{H_0}\R^3)_{\{1\}}
 \longto T^*B_0
\end{multline*}
which is a smooth Poisson manifold invariant under the induced
Hamiltonian mechanics. 
(Since $(T^*Q)_{(L_0)}|Q_{(H_1)} = \emptyset$.)
Invariance of this stratum is remarkable since this is not a general
feature of secondary strata. 

Below, when we come to explicitly describing the
induced Poisson structure we will see that this is the phase space
of a charged particle with spin moving on $B_0$. The charge is governed
by the $\R$-factor while the spin is described by the whole
$\R\times(\R^3\times\R^3))_{\{1\}}/H_0$-factor.  

The next isotropy type in the hierarchy of isotropy types of the $G$
action on $T^*Q$ is $(L_1)=(SO(2))$, and we consider $S^1=SO(2)=L_1$ as
embedded into the lower right corner of $G$. 
Now, for an element 
$x=[(k,v,w)]\in G\times_{H_0}(\R^3\times\R^3) 
\cong G\times_{H_0}\ann(\ho_0+\wo_0)$ 
the isotropy group $G_x$ with respect to the induced $G$-action
satisfies $G_x = k(H_0)_{(v,w)}k^{-1}$. Therefore,
\begin{multline*}
 \R\times(\R^3\times\R^3)_{(S^1)}/H_0\cong\R\times C^0(S^1)\hookto\\
 \hookto(T^*Q_{(H_0)})_{L_1}/G
 \cong
 (Q_{H_0}\times_{B_0}T^*B_0)\times_{W_0}(\R\times\R^3\times_{H_0}\R^3)_{(S^1)}
 \longto T^*B_0
\end{multline*}
where $C^0(S^1)\cong((\R\times\R)\setminus(0,0))/\set{\pm1}$ 
denotes the cone over $S^1$ without the cone point. This secondary
stratum is, however, not invariant under the Hamiltonian dynamics.
This means that particles with spin of $C^0(S^1)$-type may travel from
$B_0$ to $B_1$ while this is not possible for particles with
$(\R^3\times\R^3)_{\{1\}}/H_0$ spin type.

Considering the isotropy type $(L_2)=(SO(3))$ implies 
\[
 (T^*Q_{(H_0)})_{L_2}/G
 \cong
 T^*B_0\times\R\times\set{0}
\]
which, too, is not invariant under the reduced dynamics. 
The most singular element in
the isotropy lattice of $T^*Q$ is $(L_3)=(SO(4))$. However,
$(T^*Q)_{(L_3)}|_{(H_0)} = \emptyset$, and we have thus described the
secondary stratification of $(T^*Q_{(H_0)})/G$.

\subsection{Curvature, Poisson structure, and symplectic leaves}
By Section \ref{sec:curv}, in order to understand the Poisson and
symplectic structures of the reduced space $(T^*Q_{H_0})/G$ we have to
compute the mechanical connection $\curv^{\mathcal{A}}$. Thus we first
compute the inertia tensor $\ine:
Q_{H_0}=S_0^3\to\wo^*\otimes\wo^*$. That is, for $iX,iY\in i\R=\wo$ 
and $q\in S_0^3$, 
\[
 \ine_q(X,Y)
 = \vv<\dd{t}{}|_0e^{itX}q,\dd{t}{}|_0e^{itY}q>
 = -XYi^2\vv<q,q> 
 = XY
\]
which is independent of $q$.\footnote{Not knowing Proposition
  \ref{prop:curv} one would have to compute the generalized
  mechanical connection from the inertia tensor $\ine^G$ defined by
  $\ine^G_q(X,Y)=\vv<\zeta_X(q),\zeta_Y(q)>$ where $q\in Q_{(H_0)}$ and
  $X,Y\in\gu$. The resulting equations then blow up horribly.}
Thus we can identify $\R\cong i\R\cong i\R^*\cong\R^*$. 
Therefore, by diagram (\ref{d:A}),
\[
 \vv<\mathcal{A}(q,v),X>
 =
 \vv<v,\zeta_X(q)>
 =
 \vv<v,iq>X
\]
where $(q,v)\in TS^3_0$. 
That is $\mathcal{A}(q,v) = \vv<v,iq>$. To better understand this we now
use the quaternionic representation of the Hopf map 
$\eta: S^3\toto S^2$. I.e., 
$q = a+ib+j\alpha+k\beta$, $v = a'+ib'+j\alpha'+k\beta'$, and 
using the orthonormal frame
$\xi_1(q)=iq$, $\xi_2(q)=jq$, $\xi_3(q)=kq$
we
trivialize $TS_0^3 = S_0^3\times\so(3)_-$. (Here, $\so(3)_-$ denotes
the Lie algebra of right invariant vector fields on $SO(3)$.) The
frame vectors $\xi_1,\xi_2,\xi_3$ and their dual covectors
$\xi^1,\xi^2,\xi^3$ enjoy the relations
\begin{align*}
 [\xi_2,\xi_3] &= -2\xi_1,
 [\xi_3,\xi_1]  = -2\xi_2,
 [\xi_1,\xi_2]  = -2\xi_3,\\
\text{ and }
 2\xi^2\wedge\xi^3 &= d\xi^1,
 2\xi^3\wedge\xi^1  = d\xi^2,
 2\xi^1\wedge\xi^2  = d\xi^3.
\end{align*}
In terms of these we get $\mathcal{A}=\xi^1$, and therefore
$\curv^{\mathcal{A}} = d\mathcal{A} = 2\xi^2\wedge\xi^3$. Since
$\eta_0: S_0^3\toto B_0$ is a Riemannian submersion the volume form
$\nu$ -which is the standard one induced from $\R^3$- pulls back to
$\eta_0^*\nu = \xi^2\wedge\xi^3$. In other words the mechanical
curvature $\curv^{\mathcal{A}}$ drops via $\eta_0$ to
$\curv_0^{\mathcal{A}}=2\nu$. 

An immediate and easily visible consequence is the following. Assume
$\lam\in\wo^*\subset\gu^*$ is such that 
$\orb\spr{0}H_0=\set{\text{point}}$
where $\orb$ is the adjoint $G$-orbit through $\lam$. 
Then  
$(T^*Q_{H_0})/G = T^*B_0$ is a magnetic cotangent bundle equipped with
the symplectic structure $\Om^{B_0}-2\vv<\lam,\tau^*\nu>$ where $\Om^{B_0}$
is the canonical structure on $T^*B_0$
and $\tau: T^*B_0\to B_0$ is the projection.
This follows from Theorem \ref{thm:red}.

\subsection{Particular cases}\label{sub:parti}
As an interesting and representative particular case let $s\neq0$ and 
consider the
$5\times5$ matrix $\lam\in\so(5)\cong\so(5)^*$ which has 
\[
 \left(
   \begin{matrix}
    0&s\\
   -s&0
   \end{matrix}
  \right)
\]
in the top left corner and zeros elsewhere. That is, we are starting from
points of the form 
\[
 (q_0,p_0,\lam)\in (T^*Q_{(H_0)})_{L_2}/G
 \cong
 T^*B_0\times\R\times\set{0}
\]
and want to compute the singular symplectic leaves passing through
these. Therefore, by Theorem \ref{thm:red}
we are to be concerned with the singular symplectic
space $\orb\spr{0}SO(3) = \orb\cap\ho_0^{\bot}/H_0$ where 
$\orb\cong SO(5)/(S^1\times SO(3))$ 
is
the $G$-orbit through $\lam$. Now, $\lam$ can as well be written as
$\lam = s e_1\wedge e_2$ where $e_i$ is to be the standard basis vector
of $\R^5$ with $1$ in the $i$-th position and zeros elsewhere. Let
$A,B,C,D,E$ denote an arbitrary positively ordered 
orthonormal basis of $\R^5$ so that
$k=(A|B|C|D|E)$ is an arbitrary element of $SO(5)=G$. 
Then
\[
 k\lam k^{-1}
 = 
 s(ke_1)\wedge(ke_2)
 =
 s A\wedge B
\]
and the condition for $k\lam k^{-1}$ to be in $\ho_0^{\bot}$
translates to 
\[
 \left(
   \begin{matrix}
    a_3\\
    a_4\\
    a_5
   \end{matrix}
  \right)
\times
 \left(
   \begin{matrix}
    b_3\\
    b_4\\
    b_5
   \end{matrix}
  \right)
 =0,
\]
i.e.,
the vectors formed by the latter $3$ components of $A$ and $B$ are to
be linearly dependent. Via the $H_0$ action we can thus bring $A$ and
$B$ to the normal form $A=(a_1,a_2,a_3,0,0)^t$ and
$B=(b_1,b_2,b_3,0,0)^t$. 
It should be noted that this induces an action of $\set{\pm1} =
N_{H_0}(S^1)$ 
in the third
component which eliminates the remaining freedom.
In other words, the $0$-level set of the $H_0$-\momap is 
a smooth manifold
of the form
\[
 \orb\cap\ho_0^{\bot}
 \cong
 (\orb\cap\ho_0^{\bot})^{S^1}\times_{\set{\pm1}}H_0/S^1
 = 
 S^2(s)\times_{\set{\pm1}}H_0/S^1.
\]
Since the orthonormality conditions on $A$ and $B$ are preserved 
this means that $A\wedge B\in SO(3)/S^1$.  
The quotient $\orb\spr{0}H_0$ can thus be described as 
\[
 (a_1,a_2,a_3)^t\times(b_1,b_2,b_3)^t=(x_1,x_2,x_3)^t\in V(2,3)/S^1=Gr(2,3)=S^2
\]
modulo the remaining $\pm1$-action in the third component which
induces a $\pm1$-action on
$(x_1,x_2)=(a_2b_3-b_2a_3,-a_1b_3+b_1a_3)$. I.e., $\orb\spr{0}H_0$ is
the following stratified space: the regular
stratum consists of $[(x_1,x_2,x_3)^t]_{\sim}\in S_0^2(s)/\sim$ where
$(x_1,x_2,x_3)\sim(-x_1,-x_2,x_3)$ and $(x_1,x_2)\neq(0,0)$; and the
singular stratum is $\set{(0,0,s),(0,0,-s)}$.  
In fact, the singularities occurring here are only orbifold
singularities. 

This decomposition of $\orb\spr{0}H_0$ 
induces the symplectic (Sjamaar-Lerman) stratification on the
singular symplectic product space 
$T^*Q_{H_0}\times\orb\spr{0}H_0$ -- which by Theorem \ref{thm:red} can be
further reduced to 
$(T^*Q_{(H_0)})\sporb G \cong (T^*Q_{H_0}\times\orb\spr{0}H_0)\spr{0}W_0$.
Its singular (i.e., lower dimensional) 
strata $T^*Q_{H_0}\times\set{\pm s}$ 
reduce via $W$ to $T^*B_0$ with the magnetic 
symplectic form $\Om^{B_0}-2\vv<\pm s,\tau^*\nu>$. 
For the regular
stratum let 
\[
 \pi_1: T^*Q_{H_0}\times S_0^2(s)/\sim \to T^*Q_{H_0} 
 \text{ and } 
 \pi_2: T^*Q_{H_0}\times S_0^2(s)/\sim \to S_0^2(s)/\sim
\]
denote the Cartesian
projections, let $\eta: T^*Q_{H_0}\to T^*B_0$ denote the projection
associated to the horizontal lifting map with respect to $\mathcal{A}:
TQ_{H_0}\to\wo_0$, and let $\pr_0: S_0^2(s)/\sim\to\wo_0$ denote the
induced \momap of the induced Hamiltonian action by $W_0$ on $S_0^2(s)/\sim$. 
(Again, see Theorem \ref{thm:red}.)
Then the reduced symplectic form on the regular piece 
$(T^*Q_{H_0}\times\orb\spr{0}H_0)\spr{0}W_0$
is the one induced by the basic form 
\[\tag{MCF}\label{f:mcf} 
 \pi_1^*\eta^*\Om^{B_0}
 - 2\vv<\pr_0\circ\pi_2,(\pi_{W_0}\circ\tau\circ\pi_1)^*\nu>
 + \pi_2^*\gamma
\]
on $J_{W_0}^{-1}(0) \subset T^*Q_{H_0}\times S_0^2(s)/\sim$.
Here, $\gamma$ is the reduced form on the surface
$S_0^2(s)/\sim=(\orb\spr{0}H_0)_{\text{regular}}$,
and $J_{W_0} = \mu_{Q_{H_0}}-\pr_0$ is the \momap on
$T^*Q_{H_0}\times S_0^2(s)/\sim$ with respect to the diagonal
Hamiltonian $W_0$-action.
Notice that the
middle term in this formula makes the spin variables interact with the
magnetic field $2\nu$.

\subsection{Hamiltonian reduction}   
Let us carry out the reduction at level $\lam$ as in Subsection
\ref{sub:parti} once again but this time in the presence of a
Hamiltonian ${\mathcal{H}}$. The obvious $G$-invariant 
Hamiltonian to look at is the one
associated to the round metric on $S^9$, i.e.,
${\mathcal{H}}: T^*S^9\to\R$, $(q,p)\mapsto\by{1}{2}\vv<p,p>$. Use this metric to identify
$TS^9=T^*S^9$, and to split $TQ_{(H_0)}=\hor\oplus\ver$ into
horizontal and vertical parts with respect to $Q_{(H_0)}\toto
Q_{(H_0)}/G=B_0$. The generalized mechanical connection form $A$ thus
induces a collection of point-wise isomorphisms $A_q:
\ver_q\to\gu_q^{\bot}$. (See Section \ref{sec:curv}.) By appeal to the
$G$-action any point $(q,p)\in TQ_{(H_0)}$ can be transported to a point
$(q_0,p_0+\zeta_X(q_0))$ where $q_0\in Q_{H_0}$,
$p_0\in\hor_{q_0}\cong T_{[q_0]}B_0$, and
$X=\Ad_k(A_{q}(p))\in\ho_0^{\bot}$. Therefore,
\begin{align*}
 {\mathcal{H}}(q,p) 
 &= {\mathcal{H}}(q_0,p_0+\zeta_x(q_0))
  = \by{1}{2}\vv<p_0,p_0>
    + \by{1}{2}\ine_{q_0}(X,X)\\
 &= \by{1}{2}\vv<p_0,p_0>
    + \by{1}{2}\vv<X.q_0,X.q_0>\\
 &= \by{1}{2}\vv<p_0,p_0>
    + 
 \by{1}{2}\langle
 \left(
   \begin{matrix}
    s^2-x^2 & -xy & 0 & 0\\
    -xy & s^2-y^2 & 0 & 0\\
    0 & 0 & s^2-x^2 & -xy\\
    0 & 0 & -xy & s^2-y^2
   \end{matrix}
  \right)
     \left(
     \begin{matrix}
     a\\
     b\\
     \alpha\\
     \beta
     \end{matrix}
     \right),
       \left(
     \begin{matrix}
     a\\
     b\\
     \alpha\\
     \beta
     \end{matrix}
     \right)\rangle\\
 &=: 
  \by{1}{2}\vv<p_0,p_0>
    + V(q_0,X)
\end{align*}
where $X=k\lam k^{-1}\in\orb$. 
In this equation we used a result of the previous section. Namely,
that $X\in\orb\cap\ho_0^{\bot}$ can, via the $H_0$-action, be brought
to the form
\[ 
 X =
  \left(
   \begin{matrix}
    0  & z  & -y & 0 & 0\\
    -z & 0  & x  & 0 & 0\\
    y  & -x & 0  & 0 & 0\\
    0  & 0  & 0  & 0 & 0\\
    0  & 0  & 0  & 0 & 0
   \end{matrix}
  \right),
\]
and this action does not affect the value of the
Hamiltonian. Moreover, $x^2+y^2+z^2 = s^2$, and one has to distinguish
points $(x,y,z)$ according to $(x,y)=(0,0)$ on the one hand and
$(0,0)\neq(x,y)\sim(-x,-y)$ on the other hand.
Thus we have computed the reduced Hamiltonian 
\begin{multline*}
 {\mathcal{H}}_0:
 T^*Q_{(H_0)}\sporb G 
 =
 (T^*Q_{H_0}\times\orb\spr{0}H_0)\spr{0}W_0\\
 \cong
 (Q_{(H_0)}\times_{B_0}T^*B_0)\times_{W_0}\orb\spr{0}H_0\longto\R.
\end{multline*}
Since we are dealing with a stratified space we effectively get two
Hamiltonian systems. Indeed, the system corresponding to the small
stratum where $(x,y,z)=(0,0,\pm s)$ is 
\[
 \Big( 
  T^*B_0,
  \Om^{B_0}\mp s2\tau^*\nu,
  {\mathcal{H}}_0 = \by{1}{2}\vv<p_0,p_0> + V(q_0,(0,0,\pm s))
      = \by{1}{2}\vv<p_0,p_0> + \by{1}{2}s^2
 \Big),
\]
and its flow equations are the Lorentz equations 
describing the motion of a charged particle with charge $\pm s$ 
on $B_0$ under the
influence of the magnetic field $2\nu$. 
Indeed, this can be seen
as follows. 
Let $g$ denote the induced metric on $S^2(\by{1}{2})$ and $J$ the
standard complex structure such that $\hat{g}\circ J = \hat{\nu}$.
By Remark~\ref{rem:wong} (or direct computation) the Hamiltonian
equations associated to $\mathcal{H}_0$ 
(lifted to $T^*S^2(\by{1}{2})$)
yield Lorentz equations for a
curve $c(t)$ of the type 
$\nabla_{c'}c' 
 = \hat{g}^{-1}(\pm si_{c'}\curv_0^{\mathcal{A}}) 
 = \pm 2sJ(c')$, 
whence the 
curves $c(t)$ are small circles of constant geodesic curvature $\pm 2s$. 
Projecting
these to $B = B_0 \sqcup B_1$
via the $\pm 1$ action in the $z$-direction one obtains the solution
curves. In particular these solutions leave $B_0$ in finite time to be
reflected at $B_1$ back into $B_0$ in a Snell's law manner.  
                                                                                                                                                                                                                                                                                                                                                                                                                                                                                                                                                                                                                                                                                                                                                                                                                                                                                                                                                                     
Emphasizing the point
of view of minimal coupling this means that the free Hamiltonian
system on $T^*B_0$ with gauge group $W_0$ is coupled with 
with $\set{(0,0,\pm s)}\wo^*$. 
Thus the particles are equipped with
a charge $\pm s$ and accordingly deflected from their geodesic paths.    

\begin{comment}
**********************************************************************
BRAUCH ICH DAS?????????????????::
Indeed, let $h: T^*b_0\to\R$, $(q,p)\mapsto\by{1}{2}\vv<p,p>$, and
observe that the same calculation as above 
as well as the general remarks concerning Hamiltonian reduction in
Section \ref{sec:curv}
show that
$\CHI^*h = {\mathcal{H}}_0|J_{W_0}^{-1}(0)$. 
The interesting point is that by Theorem
\ref{thm:red} this coupling procedure is equivalent to breaking the
symmetries (i.e., doing
symplectic reduction) 
of a bigger, very simple system.
**********************************************************************
\end{comment}

The reduced Hamiltonian system on the regular stratum is given by
regular reduction at $0$ with respect to $W_0$ of
\[
 \Big(
  T^*Q_{H_0}\times S_0^2(s)/\sim,
  \pi_1^*\Om^{Q_{H_0}}+\pi_2^*\gamma,
  {\mathcal{H}}_0 = \by{1}{2}\vv<p_0,p_0> + V(q_0,[X]_{\sim})
 \Big)
\]
where $[X]_{\sim}\in(\orb\spr{0}H_0)_{\textup{reg}} = S_0^2(s)/\sim$. 
The corresponding reduced symplectic structure is described by 
the minimal coupling form (\ref{f:mcf}) in Subsection \ref{sub:parti}.
Let 
\[
 K[X]_{\sim}
 :=
  \left(
   \begin{matrix}
    s^2-x^2 & -xy & 0 & 0\\
    -xy & s^2-y^2 & 0 & 0\\
    0 & 0 & s^2-x^2 & -xy\\
    0 & 0 & -xy & s^2-y^2
   \end{matrix}
  \right) 
\]
and observe that 
\[
 \det K[X]_{\sim}
 = 
 (s^4-s^2(x^2+y^2))^2
 = 
 s^4z^4
 \neq0
 \iff
 z\neq0.
\]
Thus as long as $z\neq0$ the function
\[
 {\mathcal{H}}_0
 =
 \by{1}{2}\vv<p_0,p_0>
 +
 \by{1}{2}\vv<K[X]_{\sim}(a,b,\alpha,\beta)^t,(a,b,\alpha,\beta)^t>   
\] 
could be regarded as
the Hamiltonian describing a harmonic oscillator motion of the
confined charged particles $a$, $b$, $\alpha$, $\beta$ in a magnetic
field with the frequencies 
\[
 \set{\om[X]_{\sim}: \det(K[X]_{\sim}-\om[X]_{\sim})=0}
\]
being allowed some internal degrees of freedom.
(It would be interesting to know whether this is a physically relevant
example?)

\end{document}